\documentclass[a4paper,11pt]{article}
\usepackage{authblk}
\usepackage[latin1]{inputenc}
\usepackage[pdftex]{graphicx}
\usepackage{subfig}
\usepackage{float}
\usepackage{amsfonts,amssymb,amsmath,amsthm}
\usepackage{pdfsync}
\usepackage{color}
\usepackage{url}
\usepackage{nicefrac}
\usepackage[pdftex]{hyperref}
\usepackage{hyperref}

 \hypersetup{colorlinks=true,citecolor=blue, filecolor=red, linkcolor=blue, urlcolor=blue}

\newtheorem{theorem}{Theorem}[section] 
\newtheorem{proposition}[theorem]{Proposition}
\newtheorem{lemma} [theorem] {Lemma}
\newtheorem{definition}[theorem]{Definition}
\newtheorem{remark}[theorem]{Remark}

\newtheorem{hypothesis}[theorem]{Hypothesis}

\makeatletter
\@addtoreset{equation}{section}

\makeatother

\numberwithin{equation}{section}

\newcommand{\C}{\mathbb{C}}
\newcommand{\R}{\mathbb{R}}
\newcommand{\N}{\mathbb{N}}
\newcommand{\Z}{\mathbb{Z}}
\renewcommand{\l}{\ell}
\newcommand{\B}{\mathcal{B}}

\newcommand{\hexagon}{{hex}} 
\newcommand{\Zm}{\nicefrac{\Z}{6\Z}}
\renewcommand{\a}{\alpha}
\renewcommand{\b}{\beta}
\newcommand{\g}{\gamma}
\newcommand{\s}{\sigma} 
\def\ds{\displaystyle}
\def\tt{\theta}

\newcommand{\be}{\begin{equation}}
\newcommand{\ee}{\end{equation}}
\newcommand{\ben}{\begin{equation*}}
\newcommand{\een}{\end{equation*}}
\newcommand{\ba}{\begin{eqnarray}}
\newcommand{\ea}{\end{eqnarray}}
\newcommand{\ban}{\begin{eqnarray*}}
\newcommand{\ean}{\end{eqnarray*}}

 \makeindex

\begin{document}

\title{On the low lying spectrum of the magnetic Schr\"odinger operator with kagome periodicity} 
\author{Philippe Kerdelhu\'e${}^{1}$ \& Jimena Royo-Letelier${}^{2}$ }
 \affil{ \footnotesize ${}^1$   D\'epartement de Math\'ematiques, CNRS UMR 8628, \\ F-91405 Orsay Cedex, France  \\
Philippe.Kerdelhue@math.u-psud.fr\\
${}^2$ jimena.royo-letelier@m4x.org}

\date{}

\maketitle

  \begin{abstract} We study in  a semi-classical regime a periodic magnetic Schr\"odinger operator in $\R^2$.
This is inspired by recent experiments on artificial magnetism with
ultra cold atoms in optical lattices, and by the new interest for the operator on the hexagonal lattice describing the 
behavior of an electron in a graphene sheet.
We first review some results for the square (Harper), triangular and hexagonal lattices.
Then we study the case when the periodicity is given by the kagome lattice considered by Hou.
Following the techniques introduced by Helffer-Sj\"ostrand
and  Carlsson, we reduce this problem to the study  of a discrete operator on $ \ell^2(\Z^2;\C^3) $
and a pseudo-differential operator on $ L^2(\R;\C^3)$, which keep  the symmetries of the kagome lattice. We estimate the
coefficients of these operators in the case of a weak constant magnetic field.  Plotting the spectrum for rational values
of the magnetic flux divided by $2\pi h$ where $h$ is the semi-classical parameter, we obtain a picture similar to Hofstadter's butterfly. We study the properties of this
picture and prove the symmetries of the spectrum and the existence of flat bands, which do not occur in the case of
the three previous models.   \end{abstract}

\section[Introduction]{Introduction} \label{intro} 

We consider in a semi-classical regime the Schr\"odinger magnetic operator $P_{h,A,V}$, defined as the self-adjoint extension
in $L^2(\R^2)$ of the operator given in $C_0^\infty(\R^2) $ by

\be  \label{ScrhoPot3}
	P_{h,A,V}^0= (h D_{x_1}-A_1(x))^2  + (hD_{x_2}- A_2(x))^2 + V(x) \,, 
 \ee 

where $D_{x_j}=\frac1i \partial_{x_j}$.  Our goal is to study the spectrum of $P_{h,A,V}$ as a function of $A$ and
the semi-classical parameter $h>0$, when $V$ has its minima in the kagome lattice and both $V$ and $ B = \nabla \wedge A$ are
invariant by the symmetries of the kagome lattice. \\

Our interest in this  mathematical
problem is motivated by recent experiments on artificial magnetism with ultra cold atoms (\cite{DaGerJuOh,JZ}), that 
lead to new geometries for this problem. To our knowledge, the Hamiltonian in (\ref{ScrhoPot3}) has not been obtained in 
a laboratory with ultra cold atoms, but we mention that a two-dimensional kagome lattice for ultra cold atoms has been recently 
achieved (\cite{JoETall}) using optical potentials. Our  main motivation is to understand and analyze mathematically various
 considerations of Hou in \cite{Hou}. \\

Let us explain the setting of our problem. A n-dimensional Bravais lattice is the set of points spanned over $\Z$ by the 
vectors of a basis $ \{   \nu_1, \cdots , \nu_n \}$ of $\R^n$. A fundamental domain of the Bravais lattice is a domain of the
 form
\ben
	 \mathcal V  = \big\{ t_1 \nu_1 + \cdots + t_n \nu_n  \,;\,  t_1, \cdots , t_n  \in [0,1] \big\}  \,. 
 \een

\begin{figure}  [H]
	 \center \includegraphics[width=.5\textwidth]{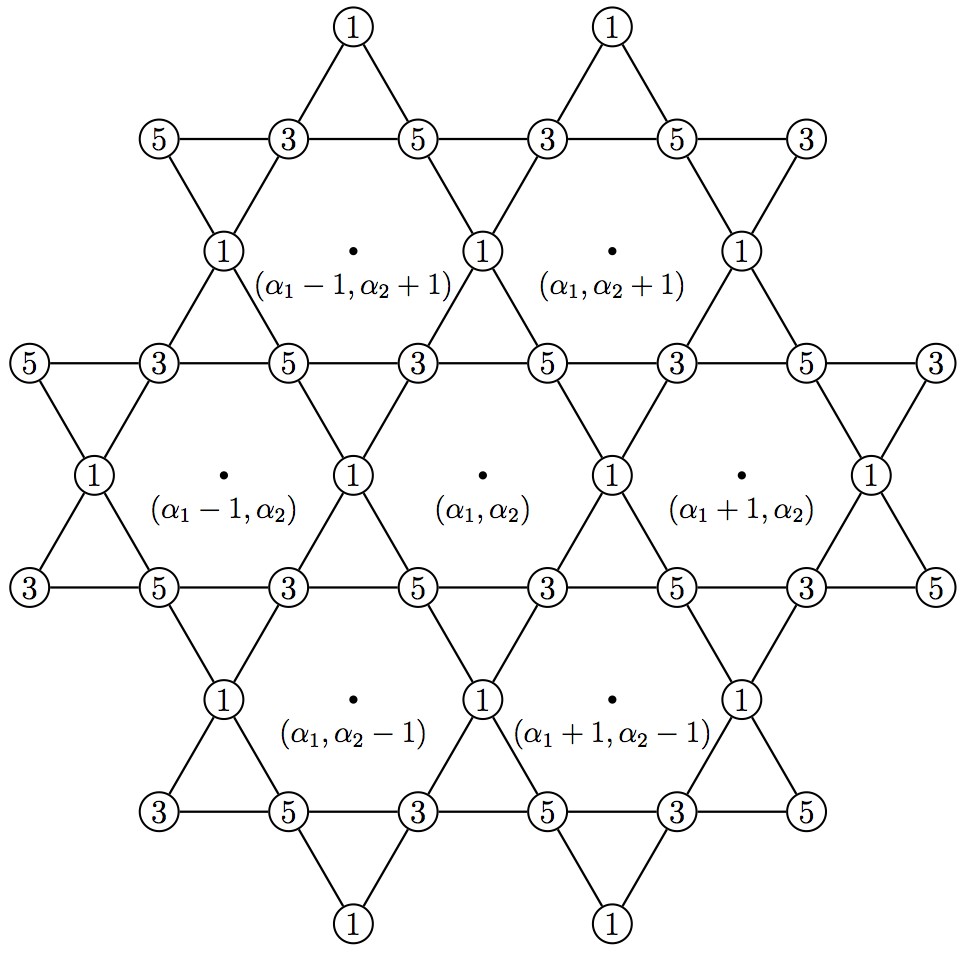}
	\caption[]{The kagome lattice and its labelling.}
\label{labeling}
\end{figure}

\begin{figure}  [H] 
  \centering
  	\subfloat[]{\label{squarelattice}\includegraphics[width= .4\textwidth]{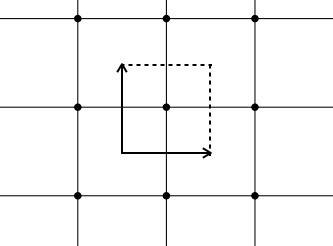}}   
  	\hspace{.5cm}
	\subfloat[]{\label{triangularlattice}\includegraphics[width= 0.4\textwidth]{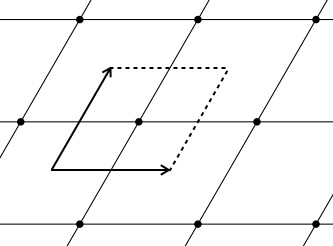}} \\
	\subfloat[]{\label{hexagonallattice}\includegraphics[width= .4\textwidth]{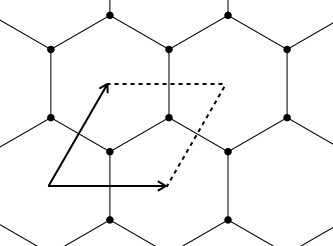}} 
	\hspace{1.5cm}
	\subfloat[]{\label{kagomelattice}\includegraphics[width=.4\textwidth]{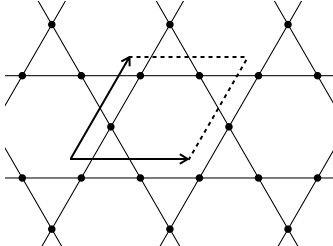}} 
	\caption[square, triangular, hexagonal and kagome lattices]{(a) square, (b) triangular, (c) hexagonal  and (d) kagome
lattices. In each case we had drawn a fundamental domain of the Bravais lattice. The points of the lattice correspond to the minima of the electric potential.}
  \label{severallattices}
\end{figure}

The kagome lattice is not a Bravais lattice, but is a discrete subset of $\mathbb R^2$ invariant under translations along a
triangular lattice and containing three points per fundamental domain of this lattice (see Figures \ref{labeling} and
\ref{kagomelattice}). Each point of the lattice has four nearest neighbours for the Euclidean distance. 
The word \textit{kagome} means a bamboo-basket (kago) woven pattern (me) and it seems that the lattice was named by the
Japanese physicist K.~Husimi in the 50's (\cite{Mekata}). \\

Let $\Gamma_{\triangle}$ be the triangular lattice spanned by $ \mathcal{B} = \{ 2 \nu_1, 2  \nu_2 \} $, where 
\be \label{defnu}
	\nu_\l = r^{\l-1}(1,0)
\ee 	

and $r$ is the rotation of angle $\nicefrac \pi3$ and center the origin. The kagome lattice can be seen as the union of 
three conveniently translated copies of $\Gamma_{\triangle}$ :  
\be  \label{defkagomelattice}
	\Gamma = \Big\{ m_{\alpha,\l} = 2 \alpha_1 \nu_1 + 2 \alpha_2 \nu_2 + \nu_\l  \,;\,  (\alpha_1,\alpha_2) \in \Z^2 \,,\, \l = 1,3,5  \Big\} \,. 
\ee  
 
We label the points of $\Gamma$ by their coordinates in $\mathcal{B}$:
\ben
	\tilde \Gamma = \Big\{ \tilde m_{\alpha,\l} = (\alpha_1, \alpha_2) + \tilde \nu_\l     \,;\, (\alpha_1,\alpha_2) \in \Z^2 \,,\, \l = 1,3,5  \Big\}     \,,
\een   

where 
\be \label{deftildenu}
	\tilde \nu_\l = \frac12 \kappa^{\l-1}(1,0)    
\ee

are the coordinates of $\nu_j$ in the basis $\mathcal{B}$. The map $\kappa : \Z^2 \to  \Z^2 $ here before is given by 
\be \label{defkappa}
	\kappa  (\alpha_1,\alpha_2  )  = (-\alpha_2,\alpha_1 +\alpha_2)  
\ee

and represents the rotation $r$ in the basis $\B$, that is, $ \widetilde{r(m_{\alpha,\l})} = \kappa(\tilde m_{\alpha,\l}) $.    \\ 

We will often consider  $j = 1,\dots,6$ as an element of  $\Zm$. Depending on the situation, we will give the properties of the kagome lattice in terms of the points $m_{\a,\l}$ or in terms of their coordinates $ \tilde m_{\a,\l}$.  \\

The symmetries of $\Gamma$ are given by those of  $\Gamma_{\triangle}$. For  $j \in \Zm $ consider the translations 
$t_j(x) = x +2\nu_j$ and define
\be \label{defG}
	\mathcal{G} = \text{the subgroup of the affine group of the plane generated by } r, \, t_1 \text{ and } t_2 \,.
\ee

Setting $(gu) (x)= u (g^{-1}(x))$  for $g \in \mathcal{G}$, we define a group action of $\mathcal{G}$ on $C^\infty(\R^2)$ which can be extended as an unitary action on $L^2(\R^2)$.    

\begin{hypothesis}    \label{hypothesisPotV}  The electric potential $V$ is a real nonnegative $ C^{\infty}$ function such that
\ba	 
	&& gV = V  \quad \text{ for all } g \in \mathcal{G}   \label{symmV} \,,  \\
	&& V \geq 0 \quad   \text{ and } \quad V(x) = 0   \text{ if and only if  } x \in\Gamma  \label{minimaV} \,, \\
	&& {\rm Hess} \, V(x) > 0 \quad \forall x \in \Gamma  \label{nodegV} \,.
\ea
 \end{hypothesis}
 
 We associated with the magnetic vector potential  $A = (A_1,A_2)$  the 1-form 
\ben
	\omega_A = A_1 dx_1 + A_2 dx_2	\,.	
\een
The magnetic field $B$ is then associated with the 2-form obtained by taking the exterior derivative of $\omega_A$:
\ben
	d  \omega_A = 	 B(x) dx_1 \wedge dx_2 \,.	
\een
In the case of $\R^2$, we identify this 2-form with $B$. The renormalized flux of $B$ through a fundamental domain
$\mathcal V$  of $\Gamma_{\triangle}$ is by definition
\ben
	\gamma  =\frac{1}{h} \int_{ \mathcal V } d\omega_A \,.
\een
 \begin{hypothesis}    \label{hypothesisPotA} 
The magnetic potential $A$ is a $C^\infty $ vector field such that the corresponding magnetic 2 form satisfies  

\be \label{symmsB}
	gB= B \quad \text{ for all } g \in \mathcal{G}   \,. 	
 \ee
\end{hypothesis}

In the case when $A=0$  (see for example Chapter XIII.16 in \cite{ReedSimons4}), the spectrum of  $P_{h,A,V}$ is continuous
and composed of bands. The general case, even when the magnetic field is constant, is very delicate. 
  The spectrum of $P_{h,A,V}$ can indeed become very singular (Cantor structure) and depends crucially on the arithmetic
properties
 of $\gamma/(2\pi)$.\\

To approach this problem, we are often led to the study of limiting models in different asymptotic regimes,
such as discrete operators defined over $\l^2(\Z^2;\C^n)$,
or equivalently, $\gamma$-pseudo-differential operators defined on $L^2(\R;\C^n)$ and associated with periodical symbols.\\ 
 
The discrete operators considered are polynomials in $\tau_1$, $\tau_2$, $\tau^*_1$ and $\tau_2^*$ with coefficients in 
$M_n(\C)$, where  $\tau_1$ and $\tau_2$ are the discrete magnetic translations  on $\ell^2(\Z^2;\C^n)$ given by 
\be \label{def1tau12}
	(\tau_1v)_\a = v_{\a_1-1,\a_2} \,, \qquad   (\tau_2v)_\a = e^{ i   \gamma \a_1} v_{\a_1,\a_2-1} \,.
\ee

We also recall that the $\gamma$-quantization of a symbol $p(x,\xi,\gamma)$ with values in $M_n(\C)$ is the pseudo-differential
operator defined over $L^2(\R ; \C^n)$ by 
\be \label{defpdo}
		\left( (\text{Op}_{\gamma}^{W}p ) u \right)(x) = \frac1{ {2\pi }\gamma} \int\!\!\! \int_{\R^2} e^{i \frac{(x-y) \xi}{ \gamma}}  p\left(\frac{x+y}2,\xi,\gamma \right) \, u(y) \, dy\,  d\xi \,.
\ee

In this article, following the ideas in \cite{HeSjHarper1}, \S 9,  we first analyze the restriction of $P_{h,A,V}$ to a spectral space
associated with the bottom of its spectrum, and we show the existence of a basis of this space such that the matrix
of this operator keeps the symmetries of $\Gamma$. \\

In order to state our first theorem, let us explain more in detail this procedure. First of all, the harmonic approximation
together with Agmon estimates shows the existence of an exponentially small (with respect to $h$) band in which one part of
the spectrum (including the bottom) of  $P_{h,A,V} $ is confined. We name this part the \textit{low lying spectrum}. The rest of the spectrum
is  separated by a gap of size $h/C$. \\

Consider  $\delta\in (0,1/8)$ and a non negative radial smooth function $\chi$, such that $\chi =1$ in $B(0,\delta/2)$
 and $\text{supp} \, \chi \subset B(0, \delta)$. For any $  m \in   \Gamma$ define 
\ben
	V_{  m} (\cdot)=  \sum_{   n \in \Gamma \setminus \{   m  \}} \chi(\cdot -n) 
\een	
and
\be \label{defPm}
	P_{m}  = P_{h,A,V}  +  V_{m} \,.
\ee 

All the $ P_{m} $ are unitary equivalent and 
\be \label{definitionb}
	b=\liminf_{|x|\to\infty}  V_{m}(x)
\ee

is positive and does not depend on $m$. The spectrum of $ P_{m} $ is discrete in the interval $[0,b]$. The first eigenvalue of $P_{m}$ is simple
and we note it $\lambda(h)$. We can prove that there exists then $\epsilon_0 >0$  such that 
$ \sigma(P_m) \cap I(h) = \{ \lambda(h) \}$, where $I(h) =    [0, h (\lambda_{har,1 }+\epsilon_0)]$ and $\lambda_{har,1 }$ is
 the first eigenvalue of the operator associated with $P_{m}$ by the harmonic approximation when $h=1$
 (see Section \ref{harappsection} for more details). We define 
\be	
	\Sigma = \text{ the spectral space associated with } I(h) \,.
\ee	

We denote by $d_V$ the Agmon distance associated with  the metric  $V dx^2$  (see \cite{DiSj}, \S 6)
 and
\be S=\min\{d_V(n,m);\, n,m\in\Gamma,\, n\not= m\}.\label{mindist} \ee 
We then have:

\begin{theorem} \label{theorem1} Under Hypotheses  \ref{hypothesisPotV} and \ref{hypothesisPotA},
there exists $h_0 > 0$ such that for $h \in (0,h_0)$ there exists a basis of $\Sigma$ in which $P_{h,A,V} \big| \Sigma$
has the matrix   
  \ben 
	\lambda(h)I + W_\gamma    \,,
  \een  
where  for all $ \tilde n, \tilde m \in \tilde  \Gamma$  and $\alpha \in \Z^2$, $W_\gamma$ satisfies 
\ba 
	(W_\gamma)_{\tilde n,\tilde m} &=&  \overline{(W_\gamma)_{\tilde m , \tilde n}}   \label{Wherm} \,,\\ 
	(W_\gamma)_{\tilde n,\tilde m} 
&=& e^{- i\frac{  \gamma}{2}(\tilde m -\tilde n) \wedge \alpha } (W_\gamma)_{(\tilde n+\alpha),(\tilde m+\alpha)}   \label{Wtrans} \,, \\
	(W_\gamma)_{\tilde n,\tilde m} &=& (W_\gamma)_{\kappa(\tilde n),\kappa(\tilde m)} \label{wk}   \,. \label{Wrot}
 \ea

 Moreover, there exists $C>0$ such that for every $\epsilon >0 $  there exists $h_\epsilon >0$, such that for
$h \in (0,h_\epsilon )$ 
 \ba \label{estcoeffintmatrix}
|(W_\gamma )_{\tilde n,\tilde m}| &\leq& C \, \exp \left( - \frac{(1-\epsilon  )  d_V (m,n)}h \right)   \label{estcoeffintmatrix} \,, \\
 |(W_\gamma )_{\tilde n,\tilde n}|  &\leq&  C \, \exp \left( - \frac{(2S-\epsilon  )}h \right)  \label{coefdiag}\,. 
  \ea
  \end{theorem}

The coefficients of $W_\gamma$ are related to the interaction between different sites of the kagome lattice. Our next result
concerns the study of this matrix, when we only keep the main terms for the Agmon distance. In order to estimate these terms,
we need additional hypothesis. Here we assume (see \cite{HeSjWsm1} for more details):  

 \begin{hypothesis}  \label{hypothesisNN} 
 \begin{itemize}
	\item[A.] The nearest neighbors for the Agmon distance are the same of those for the Euclidean distance, i.e.\, $S=d_V\left(m_{(1,0),3},m_{(0,0),1}\right)$.
	
	\item[B.] Between two nearest neighbors $m,n \in \Gamma$ there exists an unique minimal  geodesic $\zeta_{m,n}$  for  the Agmon metric.  
	
	\item[C.] This geodesic $\zeta_{m,n}$ coincides with the Euclidean one that is the segment   between $m$ and $n$. 
	
	\item[D.]  The geodesic $\zeta_{m,n}$ in non degenerate in the sense that  there is a point\footnote{Actually this condition does not depend  on the choice of the point  $  x_0 $ (see \cite{HeSjWsm1}).}   $  x_0 \in  \zeta_{m,n} \setminus \{m,n\}$ such that the function $x \mapsto d_V(x,m) + d_V(x,n) - d_V(m,n)$ restricted to a transverse line to $\zeta_{m,n}$ at $  x_0$ has a non degenerate local minimum at $  x_0$.   
\end{itemize}

 \end{hypothesis}

Under this hypothesis, we will estimate the main terms in the case of a weak and constant magnetic field $B=hB_0$, given by 
the gauge 	
\be \label{weakA}
	A(x_1,x_2) = \frac{h B_0}2 (-x_2,x_1) \,,  \quad B_0 >0 \,. 
\ee 	

The discrete model associated with the kagome lattice is 
\be \label{defQgo}
Q_{\gamma,\omega}= 
 \left(\begin{array}{ccc}
       0 & e^{i(\omega+\frac\gamma8)} \left(\tau_1^*+e^{-i\frac\gamma2}\tau_1^*\tau_2\right) & e^{-i(\omega+\frac\gamma8)} \left(\tau_1^*+\tau_2^*\right)\\
e^{-i(\omega+\frac\gamma8)} \left(\tau_1+e^{-i\frac\gamma2} \tau_1\tau_2^*\right) & 0 & e^{i(\omega+\frac\gamma8)} \left(e^{-i\frac\gamma2}\tau_1\tau_2^*+\tau_2^*\right)\\
     e^{i(\omega+\frac\gamma8)} \left(\tau_1+\tau_2\right) & e^{-i(\omega+\frac\gamma8)} \left(e^{-i\frac\gamma2}\tau_1^*\tau_2+\tau_2\right) & 0
\end{array}\right)
\ee

\normalsize
acting on $\ell^2(\mathbb Z^2;\,\mathbb C^3)$.\\

 We also introduce the symbol
\be p^{kag}(x,\xi,\gamma,\omega)= 
\left(\begin{array}{ccc}
       0 & e^{i(\omega+\frac\gamma8)} \left(e^{-ix}+e^{-i(x-\xi)}\right) & e^{-i(\omega+\frac\gamma8)} \left(e^{-ix}+e^{-i\xi}\right)\\
e^{-i(\omega+i\frac\gamma8)} \left(e^{ix}+e^{i(x-\xi)}\right) & 0 & e^{i(\omega+\frac\gamma8)} \left(e^{i(x-\xi)}+e^{-i\xi}\right)\\
     e^{i(\omega+\frac\gamma8)} \left(e^{ix}+e^{i\xi}\right) & e^{-i(\omega+\frac\gamma8)} \left(e^{-i(x-\xi)}+e^{i\xi}\right) & 0
\end{array}\right) \label{pseudodifferentialkag}
\ee

\normalsize
and its Weyl-quantization $P^{kag}_{\gamma,\omega}={\rm Op}^W_\gamma p^{kag}(x,\xi,\gamma,\omega)$  acting on $L^2(\mathbb R;\,\mathbb C^3)$.\\
 
We now state two theorems linking the Schr\"odinger operator and these two models.

\begin{theorem} \label{studyhatW}
Let $V$ satisfies Hypothesis \ref{hypothesisNN}. There exists $b_0>0$, $h_0>0$, $C>0$ and
$R\in\mathcal L(\ell^2(\mathbb Z^2;\,\mathbb C^3))$ such that for $h \in (0,h_0)$, $P_{h,A,V}|\Sigma$ is unitary equivalent with 		
\be \label{decomp1}
	Q_\gamma = \lambda(h)\,I-\rho \left( Q_{\gamma,\omega}  + R_\gamma \right)	 \,,
\ee	
where
\ba
	\rho &=&  h^{\nicefrac 12} \, b_0\, e^{-\frac Sh }  \left(1+  \mathcal{O}(h) \right) \,, \label{estrhoh} \\
	\omega &=&  \mathcal{O}(h) \label{estbetah}
\ea 
and
\be \label{estRtheo}
	 \| R_\gamma\|_{\mathcal{L}(\ell^2(\mathbb{Z}^2;\mathbb{C}^3))} \leq C \, \exp \left( - \frac{1}{Ch} \right)  \,.
\ee
\end{theorem} 

\begin{theorem} \label{theopsudodiff}
Under the same Hypothesis of Theorem \ref{studyhatW},  there exists a symbol $r(x,\xi)$ $2\pi$-periodic in $x$ and $\xi$ such
that $P_{h,A,V}|\Sigma$ has the same spectrum than 		
\be \label{decomp2}
	\lambda(h)\,I-\rho \left( P^{kag}_{\gamma,\omega}  + {\rm Op}^W_\gamma r(x,\xi)\right)	\,,
\ee	
where $\rho$ and $\omega$ are given by (\ref{estrhoh}) and (\ref{estbetah}), and 
\be
\| {\rm Op}^W_\gamma r(x,\xi)\|_{\mathcal{L}(L^2(\mathbb R;\,\mathbb{C}^3))} \leq C \, \exp \left( - \frac{1}{Ch} \right)  \,.
\ee 
\end{theorem}

 \begin{remark}
In the case of the square, triangular and hexagonal lattices, using Hypothesis 
\ref{hypothesisPotV} and \ref{hypothesisPotA}, it is possible to prove that the terms corresponding to the interaction between 
nearest neighbours of the lattices are equal, so $\omega=0$. The situation is more complex for the kagome lattice, and we are only
 able to prove equality for half of these terms. We point out that we do not see any \text{a priori} reason for equality
 between all terms, although this is assumed in some articles (\cite{Hou,HuAlt}). 

\end{remark}

We then study the dependence on $\omega$ of the spectra.

\begin{proposition} \label{gammadependence}
Let  $\sigma_{\gamma,\omega}$ be the  spectrum of $Q_{\gamma,\omega}$. We have  
\ba
\sigma_{\gamma,\omega+\frac\pi4}&=&\sigma_{\gamma-6\pi,\omega} \label{gammatrans} \,, \\
\sigma_{ \gamma ,-\omega  } &=& \sigma_{-\gamma,\omega} \,.\label{gammasym}
\ea
Thus it in enough to consider $\omega\in[0,\pi/8]$ to obtain all the spectra.
\end{proposition}

In order to compute the spectrum of $Q_{  \gamma,\omega}$, we give a last representation in the case  when $\gamma/ (2\pi)$ is a
rational number.  \\

For $p,q \in \mathbb{N}^*$ we define the matrices  $J_{p,q},K_q \in \mathcal{M}_q(\C)$ by   
 \be \label{definitpnjandk}
	J_{p,q}=\text{diag}(\exp{(2i\pi(j-1)p/q)}) \qquad \text{ and } \qquad (K_q)_{ij} =\left\{ \begin{array}{cl} 1  & \text{if }  j=i+1 \, (\text{mod } q) \\  0 & \text{if not}  
	\end{array}\right. \,.
\ee

\begin{theorem} \label{theoremmatsym}

Let $\gamma =  2\pi p/q$  with $p,q \in \N^{*}$ relatively primes and denote by $\sigma_{\gamma,\omega}$ the spectrum of $Q_{\gamma,\omega}$.
 We have   
\be 	\label{definitionsigmagamma}
 	 \sigma_{\gamma,\omega}    = \bigcup_{ \theta_1,\theta_2  \in [0,1]}  \sigma( M_{p,q,\omega,\theta_1,\theta_2}) \,, 
 \ee    

where   $ M_{p,q,\omega,\theta_1,\theta_2} \in M_{3q}(\C)$ is given by 
\be  \label{Mpqthetaijin}
 M_{p,q,\omega,\theta_1,\theta_2} = \left(\begin{array}{c|c|c} 
	&   &   \\
	  0_q   & M_{p,q,\omega,\theta_1,\theta_2}^{13}   &   M_{p,q,\omega,\theta_1,\theta_2}^{15}   \\
	  &   &   \\
	\hline 
	  &   &   \\
	 \left[  M_{p,q,\omega,\theta_1,\theta_2}^{13} \right]^*   &   0_q   & M_{p,q,\omega,\theta_1,\theta_2}^{35}    \\
	  &   &   \\
	\hline 
	 &   &   \\
	\left[  M_{p,q,\omega,\theta_1,\theta_2}^{15} \right]^* & \left[  M_{p,q,\omega,\theta_1,\theta_2}^{35} \right]^*  &   0_q   \\
	  &   &    
	\end{array}\right)
\ee 

with
 \ba 
M_{p,q,\omega,\theta_1,\theta_2}^{13}   &=& e^{i( \omega+\frac{\pi}{4}\frac pq)} \left(  e^{i2\pi\theta_1} K_q +e^{-i\pi\frac pq} e^{i  2\pi (\theta_1+\theta_2)} K_q  J_{p,q}   \right) \nonumber\\
M_{p,q,\omega,\theta_1,\theta_2}^{15}   &=& e^{-i( \omega+\frac{\pi}{4}\frac pq)} \left( e^{i2\pi\theta_1} K_q  +e^{-i2\pi\theta_2} J^*_{p,q}  \right)  \label{Mpqtheta12in} \\
M_{p,q,\omega,\theta_1,\theta_2}^{35}   &=
& e^{i( \omega+\frac{\pi}{4})\frac pq} \left( e^{- i\pi \frac pq } \, e^{ -2i\pi (\theta_1+\theta_2)} K_q^* J^*_{p,q} + e^{- i2\pi\theta_2} J^*_{p_q}  \right) \,.     \nonumber
\ea 
\end{theorem}

\begin{remark} Formally we obtain (\ref{pseudodifferentialkag}) and (\ref{Mpqthetaijin}) by replacing the pair of
 operators $(\tau_1,\tau_2)$ in (\ref{defQgo}) by $({\rm op}^W_\gamma (e^{ix}),{\rm op}^W_\gamma (e^{i\xi}))$ and 
$(e^{-i2\pi\theta_1}K_q^*,e^{i2\pi\theta_2}J_{p,q})$.
Note that these pairs of operators have the same commutation relation
\ban \tau_2\tau_1&=&e^{i\gamma} \tau_1\tau_2 \,, \\
{\rm op}^W_\gamma(e^{i\xi}) \,
{\rm op}^W_\gamma (e^{ix}) &=& e^{i\gamma} \, {\rm op}^W_\gamma(e^{ix})\, {\rm op}^W_\gamma(e^{i\xi}) \,, \\
(e^{i2\pi\theta_2}J_{p,q})( e^{-i2\pi\theta_1}K_q^*) &=& e^{i\gamma} (e^{-i2\pi\theta_1}K_q^* )(e^{i2\pi\theta_2}J_{p,q}) \,,\ean
and we obtain three isospectral operators $Q_{\gamma,\omega}$, $P^{kag}_{\gamma,\omega}$ and $\mathcal M_{p,q,\omega}$
where $\mathcal M_{p,q,\omega}$ acts on $L^2([0,1]^2;\mathbb C^{3q})$ by 
\ben (\mathcal M_{p,q,\omega} u)(\theta_1,\theta_2)=M_{p,q,\omega,\theta_1,\theta_2}u(\theta_1,\theta_2)\,. \een
In the formalism of rotational algebras, it is said that these three isospectral operators
are representations of the same Hamiltonian in different rotation algebras (see \cite{BellissardKreftSeiler}).
\end{remark}

\begin{figure} [H]
\centering
	 \includegraphics[width=1\textwidth]{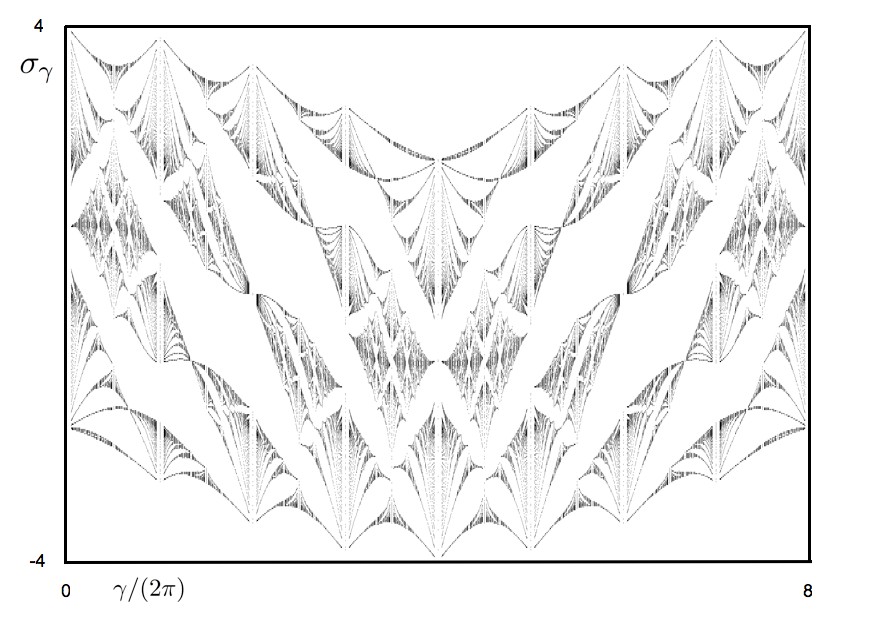} 
	\caption[Hofstadter's butterfly for the kagome lattice]{Hofstadter's butterfly for the kagome lattice when $\omega=0$.}
	\label{spectKagome}
\end{figure}

\begin{figure} [H]
\centering
	 \includegraphics[width=1\textwidth]{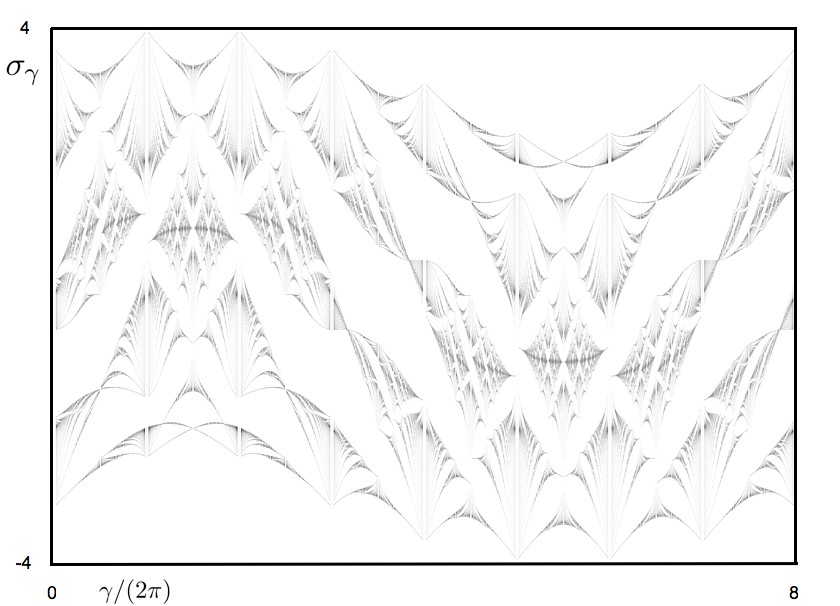} 
	\caption[Hofstadter's butterfly for the kagome lattice]{Hofstadter's butterfly for the kagome lattice when $\omega=\pi/8$.}
	\label{spectKagome8}
\end{figure}

In Figures \ref{spectKagome} and \ref{spectKagome8} we present the  equivalent of  Hofstadter's butterfly for the kagome lattice
 in the case when $\omega=0$ and $\omega=\pi/8$, obtained by numerically diagonalizing  the matrices
$ {M}_{p,q,0,\theta_1,\theta_2}$ and $ {M}_{p,q,\frac{\pi}{8},\theta_1,\theta_2}$. In the first case we recover that 
one obtained by Hou in \cite{Hou}.\\

We notice that for fixed $\gamma= 2\pi p/q$ the spectrum is composed of $3q$ (possibly not disjoint) bands,
which are the images of
\be 
	 [0,1] \times [0,1] \ni(\theta_1,\theta_2) \mapsto  \lambda_{p,q,\omega,\theta_1,\theta_2}^k \,, \quad  1\leq k \leq 3q 
\ee
where $ \lambda_{p,q,\omega,\theta_1,\theta_2}^k$ is the $k$th eigenvalue of $M_{p,q,\omega,\theta_1,\theta_2}$.\\
Since the smallest positive integer for which the operator $Q_{\gamma,\omega}$ is  invariant by the transformation
$\gamma \mapsto \gamma+2\pi k$ is $k=8$, we plot in the vertical axis of Figures \ref{spectKagome} and \ref{spectKagome8}
the bands of the spectrum  for 
\ben
	\frac{\gamma}{ 2\pi}  =  \frac pq   \,, \quad  p,q  \text{ relatively  prime and } 0\leq p < 8q \leq 400\,.
\een
   
We first observe some symmetries in these butterflies and prove the proposition 

\begin{proposition} \label{propsymmkag}

Let  $\sigma{_\gamma,\omega}$ be the  spectrum of $Q_{\gamma,\omega}$. We have  
\ba
	 \sigma_{\gamma,\omega }  & \subset & [-4,4] \label{kagrange} \,,\\
	\sigma_{\gamma +16 \pi,\omega} &=& \sigma_{\gamma,\omega}  \hspace{1.25cm} (\text{translation invariance}),  \label{kagtrans} \\
	e \in \sigma_{ \gamma + 8\pi,\omega  } &\Leftrightarrow& -e \in  \sigma_{\gamma,\omega} \quad  (\text{translation anti-invariance}).  \label{kagreflex4} 
%
\ea
In the case when $\omega=0$, we have
\ba
\sigma_{-\gamma,0} &=& \sigma_{\gamma,0 }  \hspace{1.25cm} (\text{reflexion with respect to the axis $\gamma=0$}), \label{kagreflex1} \\
e \in \sigma_{ 8\pi  - \gamma,0  } &\Leftrightarrow& -e \in  \sigma_{\gamma,0} \quad  (\text{reflexion with respect  to the
point $ (4\pi,0)$}). \label{kagreflex5}
\ea
In the case when $\omega=\pi/8$, we have
\ba
\sigma_{ 6\pi  - \gamma,\frac\pi8  } &=& \sigma_{\gamma ,\frac\pi8  } \hspace{1.25cm} (\text{reflexion with respect  to the axis $ \gamma=3\pi$}), \label{kagreflex11} \\
e \in \sigma_{-2\pi-\gamma,\frac\pi8} &\Leftrightarrow&-e \in  \sigma_{\gamma ,\frac\pi8} \quad (\text{reflexion with respect
to the point $(-\pi,0)$}). \label{kagreflex12} 
\ea
\end{proposition}

\begin{figure}   [H]
 \centering
  	\subfloat[]{\label{bands011}\includegraphics[width=.45\textwidth]{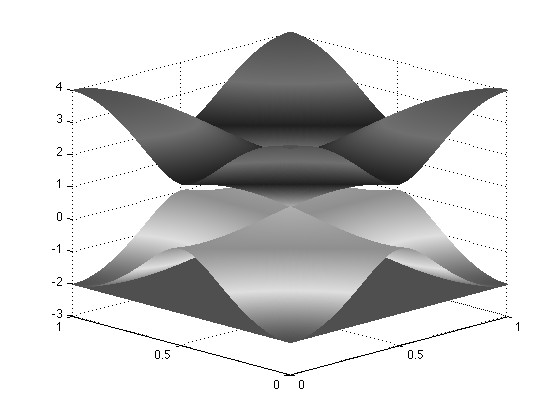} }
	\hspace{.75cm}
     	\subfloat[]{\label{bands211}\includegraphics[width=.45\textwidth]{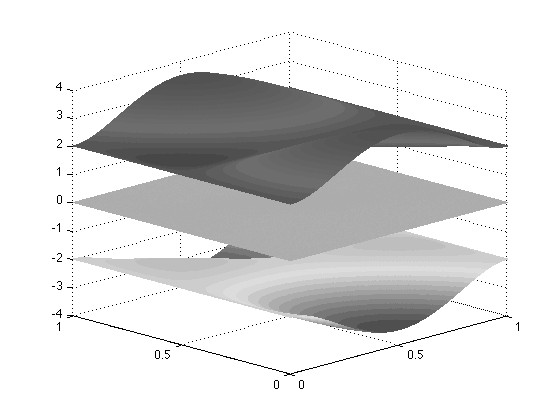} }
	\caption{Spectrum bands of $Q_{\gamma,\omega}$ for (a) $(\gamma,\omega)=(0,0)$ and (b) $(\gamma,\omega)=(4\pi,0)$.}
 	\label{bandskagome1}
\end{figure}

\begin{figure}    [H]
 \centering
  	\subfloat[]{\label{bands01}\includegraphics[width=.45\textwidth]{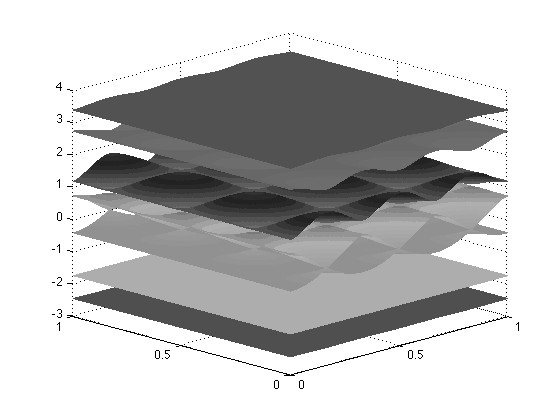} }
	\hspace{.75cm}
     	\subfloat[]{\label{bands21}\includegraphics[width=.45\textwidth]{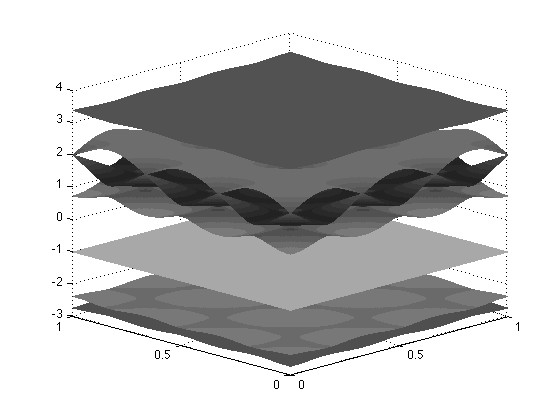} }
	\caption{Spectrum bands of $Q_{\gamma,\omega}$ for (a) $(\gamma,\omega)=(4\pi/3,0)$ and (b) $(\gamma,\omega)=(8\pi/3,0)$.}
 	\label{bandskagome}
\end{figure}

Second we note the presence of isolated points in $\sigma_{\frac{4\pi}{3},0}$,
$\sigma_{\frac{8\pi}{3},0}$,  $\sigma_{4\pi,0}$,
$\sigma_{\pi,\frac\pi8}$, $\sigma_{3\pi,\frac\pi8}$,
$\sigma_{\frac{7\pi}3,\frac{\pi}{8}}$ and $\sigma_{-\frac{\pi}{3},\frac {\pi}{8}}$.\\
To see more precisely
 the last phenomenon, we plot in Figures \ref{bands211} and \ref{bandskagome} the bands of the spectra
$\sigma_{4\pi,0}$, $\sigma_{\frac{4\pi}{3},0}$ and $\sigma_{\frac{8\pi}{3},0}$.

Numerically it seems that the second, third and fourth bands of 
$\sigma_{\frac{4\pi}{3},0}$ are reduced to $\{-\sqrt{3}\}$, and that the third, fourth and fifth band of
$\sigma_{\frac{8\pi}{3},0}$ are reduced to $\{-1\}$.\\
This leads to the definition

\begin{definition} \label{defflat}
 Let $\lambda_0$ be a reel number and $n$ a positive integer. $\{\lambda_0\}$ is called a flat band of multiplicity $n$ of 
$\sigma_{\gamma,\omega}$ if the $k^{th}$ band of $\sigma_{\gamma,\omega}$ is reduced to $\{\lambda_0\}$ for exactly $n$
values of $k$.
\end{definition}

One can easily compute 
the characteristic polynomials of the $3\times 3$ matrices $M_{0,1,0,\theta_1,\theta_2}$ and $M_{2,1,0,\theta_1,\theta_2}$.
For the other cases, we use the symbolic computation software Mathematica and
obtain

\begin{proposition} \label{eigenvalue}
 \begin{enumerate}
  \item 
\begin{enumerate}
\item $\{-2\}$ and $\{0\}$ are flat bands of multiplicity $1$ of $\sigma_{0,0}$ and
$\sigma_{4\pi,0}$ respectively.
$\sigma_{0,0}$ is composed of the three touching bands
$\{-2\}$, $[-2,1]$ and $[1,4]$.
$\sigma_{4\pi,0}$ is composed of the three disjoint bands
$[-2\sqrt{3},-\sqrt{3}]$, $\{0\}$ and $[\sqrt{3},2\sqrt{3}]$.
\item $\{-\sqrt{3}\}$ and $\{-1\}$ are flat bands of multiplicity $3$ of $\sigma_{\frac{4\pi}{3},0}$ and 
 $\sigma_{\frac{8\pi}{3},0}$ respectively.
\end{enumerate}
\item \begin{enumerate}
\item $\{-\sqrt{2}\}$ and $\{-2\}$ are flat bands of multiplicity $2$ of $\sigma_{\pi,\frac{\pi}{8}}$ and
 $\sigma_{3\pi,\frac {\pi}{8}}$ respectively. $\sigma_{3\pi,\frac {\pi}{8}}$ is composed
of the flat band $\{-2\}$ and the four touching bands\\
$[1-\sqrt{6},1-\sqrt{3}]$, $[1-\sqrt{3},1]$,
$[1,1+\sqrt{3}]$, and $[1+\sqrt{3},1+\sqrt{6}]$.
\item $\{ -\frac{\sqrt{6}+\sqrt{2}}2\}$ and $\{-\frac{\sqrt{6}-\sqrt{2}}2\}$ are flat bands of multiplicity $6$
of $\sigma_{\frac{7\pi}3,\frac{\pi}{8}}$ and $\sigma_{-\frac{ \pi}{3},\frac {\pi}{8}}$ respectively.
 \end{enumerate}
\end{enumerate}
\end{proposition}

\begin{remark} This phenomenon does not occur for the square, triangular and hexagonal models.
\end{remark}

\begin{remark} 1. Proposition \ref{eigenvalue} ensures the existence of eigenvalues of infinite multiplicity for 
$Q_{\gamma,\omega}$ and $P^{kag}_{\gamma,\omega}$ for several values of $(\gamma,\omega)$.\\
2. Since the models $Q_{\gamma,\omega}$ and $P^{kag}_{\gamma,\omega}$ only take into
 account the interactions beetween nearest wells, and $\omega=\mathcal{O}(h)$ does not a priori vanish, the existence of 
eigenvalues for $Q_{\gamma,0}$ when $\gamma$ equals to $4\pi/3$, $8\pi/3$ or $4\pi$ does not imply the existence of eigenvalues
 for the corresponding initial Schr\"odinger
operator $P_{h,A,V}$. However, Proposition \ref{eigenvalue} together with Theorem \ref{studyhatW} ensure that, when the values
of $A$ and $h$ lead to one of these values of $\gamma$, there exists $C>0$ suth that  
 a part of the low lying spectrum of $P_{h,A,V}$ is included in an interval of length at most
$C\,h^{\nicefrac32}\exp(-S/h)$ and separated from the rest of the spectrum by intervals of lengh at least
$C^{-1}\,h^{\nicefrac12}\exp(-S/h)$. 
\end{remark}

\begin{remark}
 In the light of Proposition \ref{eigenvalue} we can state the following conjecture : if $\sigma_{2\pi p/q,\omega}$
contains a flat band for a real number
$\omega$ and two relatively prime integers $p$ and $q$ with $q>0$, then 
its multiplicity is $q$.
\end{remark}

An interesting question is to see how the invariances of the initial problem are conserved in the reduced model $p^{kag}$.
The invariance by rotation of angle $\pi/3$ gave the application $\kappa$ on the indices $\alpha$, so the transpose application
${}^t\!\kappa(x,\xi)=(\xi,-x+\xi)$ is seen as the rotation of angle $-\pi/3$ on the phase space $\mathbb R_x \times \mathbb R_\xi$.
We introduce the translations $\tilde t_1(x,\xi)=(x+2\pi)$
and $\tilde t_2(x,\xi)=(x,\xi+2\pi)$, and the symmetry $s(x,\xi)=(\xi,x)$. We then have 

\begin{proposition} \label{symmsymb}
\ba p^{kag}\circ \tilde t_1=p^{kag} \,, &&  \\
p^{kag}\circ \tilde t_2=p^{kag}  \,, &&\\
\left(\begin{array}{ccc}
 0&1&0\\
0&0&1\\
1&0&0
\end{array}\right)^{-1}
\left( p^{kag}\circ {}^t\!\kappa^{2}\right)
\left(\begin{array}{ccc}
 0&1&0\\
0&0&1\\
1&0&0
\end{array}\right)&=&p^{kag} \,, \\
\left(\begin{array}{ccc}
 0&0&1\\
0&1&0\\
1&0&0
\end{array}\right)
\overline{\left( p^{kag}\circ s\right)}
\left(\begin{array}{ccc}
 0&0&1\\
0&1&0\\
1&0&0
\end{array}\right)
&=&p^{kag}\,.
\ea
\end{proposition}




\begin{remark} 1. The invariance by the rotation of angle $\pi/3$ seems lost, but fortunately the action of
 the group generated by $\tilde t_1$, $\tilde t_2$, $s$ and ${}^t\!\kappa^2$ on the set of the microlocal wells,
which are at energy $\lambda$ the connected components of $\{(x,\xi);\,\det(\lambda\,I_3 - p^{kag}(x,\xi,0,\omega))=0\}$,
is transitive.\\
2. As shown in \cite{kerdelhue2}, the invariances of $p^{kag}$ give operators commuting with
$P_{\gamma,\omega}^{kag}$.\\
We will develop these points in a further work joint with B.Helffer and devoted to the microlocal study of
$P^{kag}_{\gamma,\omega}$.
\end{remark}

\paragraph{Outline of the article }~\\

This article is organized as follows: 

\begin{itemize}
\item In Section \ref{persymb} we present a general theorem on the Weyl quantization of periodical symbols.	

\item In Section \ref{review} we review the cases of the square, triangular and hexagonal lattices. We describe and prove the
 symmetries of the corresponding spectra. 
		
	\item  In Section \ref{lattice} we study the properties of the kagome lattice and we construct a family of potentials
 invariant by the symmetries of  $\Gamma_{\triangle}$, whose minima are located in  $\Gamma$. 
	
	\item Section \ref{operator} is devoted to the semi-classical analysis  of the low lying spectrum of $P_{h,A,V}$ 
for $h$ small. We derive the discrete operator $W_\gamma $ and prove Theorem \ref{theorem1}. 
	
	\item In Section \ref{matrix}, we study the properties of $  W_\gamma $ and prove Theorem 
\ref{studyhatW}.  We give the representation using the pseudo-differential operator acting on $ L^2(\R;\C^3) $ and prove 
Theorem \ref{theopsudodiff}.  We then study the case when  $ \gamma/(2\pi)$ is rational and prove Theorem \ref{theoremmatsym}.
 We end this article by proving Propositions \ref{gammadependence}, \ref{propsymmkag} and \ref{eigenvalue}.

\end{itemize}

\noindent{\bf Ackowledgements :} This article is a revisited version of the second part of J. Royo-Letelier's PhD thesis
(defended in June 2013 at the Universit\'e de Versailles Saint-Quentin-en-Yvelines) with B. Helffer as advisor and written
with the help of P. Kerdelhu\'e.  We warmly thank B. Helffer for suggesting us this problem and for his precious help with the
realization of this article.
The second author thanks the Institute of Science and Technology Austria (IST Austria) in which she was staying
as a post-doc while this article was finalized.
The first author thanks P. Gamblin for useful conversations.

\section{Quantization of periodical symbols} \label{persymb}

We first give a general theorem on the $\gamma$-quantization of a periodic symbol, which will be used to study the symmetries of
the butterflies associated to the square, triangular, hexagonal and kagome models, and in the proofs of Theorems
\ref{studyhatW} and \ref{theopsudodiff}.
This theorem was first established in \cite{HeSjHarper1} and
\cite{kerdelhue1} for Harper's and triangular models, and under the
restriction $0<\gamma<2\pi$. We present here a slightly different proof to avoid this restriction.  \\

Let $n\in\N^*$ and $(\b,\g)\mapsto p_{\b,\g}$ be a function on $\Z^2\times\R^*$ with values in $M_n(\C)$ such that :
\ba
	\forall N\in\N,~~\exists C_N>0,~~\forall (\b,\g)\in\Z^2\times\R^*, && |p_{\b,\g}|\leq C_N (1+\b_1+\b_2)^{-N}  \label{propsymbol1} \,, \\
	\forall (\b,\g)\in\Z^2\times\R^*, && p_{-\b,\g}=p^*_{\b,\g} \,.  \label{propsymbol2}
\ea

We define the symbol
\be \label{defsymbolp}
	p(x,\xi,\g)=\sum_{\b\in\Z^2} p_{\b,\g}\,e^{i(\b_1x+\b_2\xi)} \
\ee 	
and its Weyl quantization $P_{\g}$ introduced in (\ref{defpdo}). A straightforward computation gives that $P_{\g}$ acts on
$L^2(\mathbb R;\,\mathbb C^n)$ by  
\be \label{defGgamma}
	P_{\g}\,u(x)= \sum_{\b\in\Z^2} p_{\b,\g}\,e^{i\frac{\g}{2}\b_1\b_2}\,e^{i\b_1x}\,u(x+\b_2\g) \,.
 \ee 
  
We also consider the discrete operator
\ben Q_\gamma=\sum_{\b\in\Z^2} p_{\b,\g}\,e^{i\frac{\g}{2}\b_1\b_2}\,\tau_1^{\b_1}\tau_2^{\b_2} \een

where $\tau_1$ and $\tau_2$ are the discrete magnetic translations defined in (\ref{def1tau12}), and $A_{\g}$
the infinite matrix defined by  
\ben
	\left(A_{\g}\right)_{\a,\b}=e^{-i\frac{\g}{2}\a\wedge\b}p_{\a-\b,\g} 
\een

and acting on $\ell^2(\Z^2;\,\C^n)$ by
\ben
	\left(A_\g v\right)_\a=\sum_{\b\in\Z^2}  \left(A_\g\right)_{\a,\b} \, v_\b  \,.
\een

\begin{theorem}  \label{theopseudomatrix}
	
$A_{\g}$ and $Q_\gamma$ are unitary equivalent. $P_{\g}$, $A_{\g}$ and $Q_\gamma$ have the same spectrum.

\end{theorem}
 
 {\bf Proof.} The first hypothesis anables to prove the convergence of the series defining $p$, $A_\gamma$ and $Q_\gamma$,
and the second one gives the self-adjointness of $A_\gamma$, $Q_\gamma$ and $P_{\g}$.\\

$Q_\gamma$ acts on $\ell^2(\mathbb Z:\mathbb C^n)$ by
\ban Q_\gamma u(\alpha)
&=&\sum_{\b\in\Z^2} p_{\b,\g}\,e^{-i\frac{\g}{2}\b_1\b_2} e^{i\gamma \a_1\b_2} u_{\a-\b}\\
&=&\sum_{\b\in\Z^2} p_{\a-\b,\g}\,e^{i\frac{\g}{2}(\a_1\a_2-\b_1\b_2-\a\wedge\b)} u_\b
\ean
so $A_{\g}$ and $Q_\gamma$ are unitary equivalent.\\

The operator $P_{\g}$ commutes with the translation $u(\cdot)\mapsto u(\cdot-2\pi)$, so Floquet theory applies
and the spectrum of $P_{\g}$
is the union over $\tt\in\R$ of the spectra of the operators $P^{\tt}_{\g}$ acting on the space
$\ds \left\{u\in L_{\text{loc}}^2(\R;\C^n)\,;\,u(\cdot+2\pi)=e^{i2\pi\tt}u(\cdot)\,a.e. \right\}$ by
\ben
	P^{\tt}_{\g}\,u(x)=\sum_{\b\in\Z^2} p_{\b,\g}\,e^{i\frac{\g}{2}\b_1\b_2}\,e^{i\b_1x}\,u(x+\b_2\g) \,.
\een
	
We notice that  $P^{\tt}_{\g}$ has the same spectrum than its conjugate $\tilde P^{\tt}_{\g}=e^{-i\tt x}P^{\tt}_{\g}\,e^{i\tt x}$
acting on $L^2(\R/2\pi\Z; \C^n)$ by
\ben
	\tilde P^{\tt}_{\g}\,u(x)=\sum_{\b\in\Z^2} p_{\b,\g}\,e^{i\frac{\g}{2}\b_1\b_2}\,e^{i\b_1x}\,e^{i\g \b_2\tt}u(x+\b_2\g) \,.
\een
	
The union over $\tt\in\R$ of the spectra of the operators $\tilde P^{\tt}_{\g}$ is the union over $\tt\in[0,2\pi/\g ]$
(or $\tt\in[2\pi/\g ,0]$ in the case when $\g<0$) of these spectra. Hence the spectrum of $P_{\g}$ is the spectrum of the operator
$\check P_{\g}$
acting on
$ L^2(\R/2\pi\Z\times \R/\frac{2\pi}{\g}\Z; \C^n)$ by

$$\check P_{\g}u(x,\tt)=\sum_{\b\in\Z^2} p_{\b,\g}\,e^{i\frac{\g}{2}\b_1\b_2}\,e^{i\b_1x}\,e^{i\g \b_2\tt}u(x+\b_2\g,\tt) \,. $$

We define the unitary Fourier transform $\mathcal F$ mapping $ L^2((\R/2\pi\Z)\times (\R/\frac{2\pi}{\g}\Z); \C^n)$
on $\ell^2(\mathbb Z^2;\mathbb C^n)$ by
\ben (\mathcal F u)_\alpha
=\frac{\gamma^{\nicefrac12}}{2\pi} \int\!\!\!\int e^{-i(\alpha_1 x+\gamma\alpha_2\theta)} u(x,\theta)\,dx\,d\theta \een
and a straightforward computation gives 
\ben \mathcal F\,\check P_\gamma=Q_\gamma \,\mathcal F \,. \een

So  $ \check P_{\g}$ and $Q_{\g}$ are unitary equivalent. Hence  $P_{\g}$ and $Q_{\g}$ have the same spectrum.\\  \qed\\

\section {The square, triangular and hexagonal lattices} \label{review}
  \subsection{Presentation of the models}

The spectral properties of $P_{h,A,V}$ have been studied for the square, triangular and hexagonal lattices. When plotting the
spectrum as a function of $\gamma$, we obtain a picture with several symmetries, which are determined by the symmetries of the
lattice. In the case of the square lattice, we get the famous Hofstadter butterfly. In this section we review and prove the
symmetries of these spectra using the pseudo-differential operators associated with these lattices. We recall that the symbols 
corresponding to the square, triangular and hexagonal lattices are respectively
\ba  
	p^\Box(x,\xi) &=& \cos x + \cos \xi \,, \label{defpsudosqr}  \\  
	p^\triangle(x,\xi) &=& \cos x + \cos \xi +\cos (x-\xi) \,, \label{triangularsymbol} \\ 
	p^\hexagon(x,\xi) &=&  \left(\begin{array}{cc} 0  & 1 + e^{ix} + e^{i\xi}\\ 1 + e^{-ix} + e^{-i\xi} & 0\end{array}\right) \,. \label{hexagonalsymbol}
 \ea 
  
The size of the matrix is the number of points of the lattice in each fundamental domain. The number of terms of the form $e^{i(ax+b\xi)}$ is the product of this number by the number of nearest neighbours of each point of the lattice.\\

To study the symmetries of the Hofstadter's butterflies associated with each model, we will
use the following direct consequence of Theorem \ref{theopseudomatrix}.

\begin{proposition} \label{propositionsym}

Let $n\in\N^*$ and $\beta\mapsto p_{\b}$ be a function on $\Z^2$ (here $p_{\b}$ does not depend on $\g$) with values in $M_n(\C)$ satisfying (\ref{propsymbol1}) and
(\ref{propsymbol2}). Consider the symbol $ p(x,\xi)$ defined in (\ref{defsymbolp}), its Weyl quantization $P_{\g}$ and denote by
$\s(P_{\g})$ the spectrum of $P_{\g}$.  We have:
\ba
	 \text{ If }  p_{\b}=0 \text{ for  } \b_1 \text{ and } \b_2 \text{ odd, then \, } && \forall \g\in\R \,,  \s(P_{\g+2\pi})=\s(P_{\g}) \,. \label{propodd} \\
	 \text{ If }  p_{\b}=0 \text{ for  } \b_1 \text{ and } \b_2 \text{ even, then } && \forall \g\in\R \,,  \s(P_{\g+2\pi})=-\s(P_{\g}) \,. \label{propeven}
\ea
\end{proposition}

\begin{remark} Property (\ref{propodd}) applies to the Harper model $p^\Box$ and the hexagonal model $\ds p^\hexagon $.
Property (\ref{propeven}) applies to the triangular model $p^\triangle$.
\end{remark}

{\bf Proof.} First we notice that the magnetic translations $\tau_1$ and $\tau_2$ defined in (\ref{def1tau12}) don't change
when we replace $\gamma$ by $\gamma+2\pi$. Hence Theorem \ref{theopseudomatrix} gives 

$$\s(P_{\g+2\pi})= \s\left( {\rm Op}^W_{\g}\left(\sum_{\b\in\Z^2} (-1)^{\b_1\b_2}p_{\b}\,e^{i(\b_1x+\b_2\xi)}\right)\right) \,,$$

so  (\ref{propodd}) is proved.\\

Since  the application $(x,\xi)\mapsto(x+\pi,\xi+\pi)$ is affine and symplectic the operators  
$$   {\rm Op}_{\g}^W\left(\sum_{\b\in\Z^2} (-1)^{\b_1\b_2}p_{\b}\,e^{i(\b_1x+\b_2\xi)}\right) \quad \text{ and } \quad      {\rm Op}_{\g}^W\left(\sum_{\b\in\Z^2} (-1)^{\b_1\b_2}p_{\b}\,e^{i(\b_1(x+\pi)+\b_2(\xi+\pi))}\right)$$
are unitary equivalent. Then, 
\ban
\sigma\left(\sum_{\b\in\Z^2} (-1)^{\b_1\b_2}p_{\b}\,e^{i(\b_1x+\b_2\xi)}\right)
&=&\s\left(\sum_{\b\in\Z^2} (-1)^{\b_1\b_2}p_{\b}\,e^{i(\b_1(x+\pi)+\b_2(\xi+\pi))}\right)\\
&=&-\s\left(\sum_{\b\in\Z^2} (-1)^{(\b_1+1)(\b_2+1)}p_{\b}\,e^{i(\b_1x+\b_2\xi)}\right) \,, \ean
 which yields (\ref{propeven}). \\ \qed

\subsection{The square lattice}

The square lattice is the Bravais lattice associated with the basis $\{(1,0), (0,1) \}$ of $\R^2$. Each point of the lattice
has 4 nearest neighbours for the Euclidean distance. One of the models used in this case is the discrete operator
$L^\Box_{ \gamma} $ defined on $\ell^2(\Z^2,\C)$ by     
\be \label{opharperC5}
	L^\Box_\gamma =  \frac 12 \left( \tau_1 + \tau_1^* +  \tau_2 + \tau_2^* \right)  \,,
\ee 

where $\tau_1,\tau_2$ are the discrete magnetic translations  defined in (\ref{def1tau12}).\\

Using a partial Floquet theory\footnote{The classical reference for Floquet theory is \cite{ReedSimons4}, \S XII.16. We also refer
to the review about periodic operators in Subsections 2.1 and 2.2 of \cite{PSTd}.}, we are led to the study of  the spectrum of a
family (parametrized by $\theta_2$) of discrete Schr\"odinger operators  $L^\Box_{ \gamma, {\theta_2}} $ acting over 
$\l^2(\Z)$ by
\be \label{opharperC5red}
	(L^\Box_{ \gamma, {\theta_2}} v)_n = \frac{v _{n+1} +  v _{n-1}}2  + V_{{\theta_2}} (n)  v_n \,,
 \ee 

where $V_{\theta_2} (n)=  \cos \left( \gamma n + \theta_2 \right)  $ is the discrete potential. \\   

\begin{figure}  [H]
  \centering
  	\subfloat[]{\label{spectSquare}\includegraphics[width=.5\textwidth]{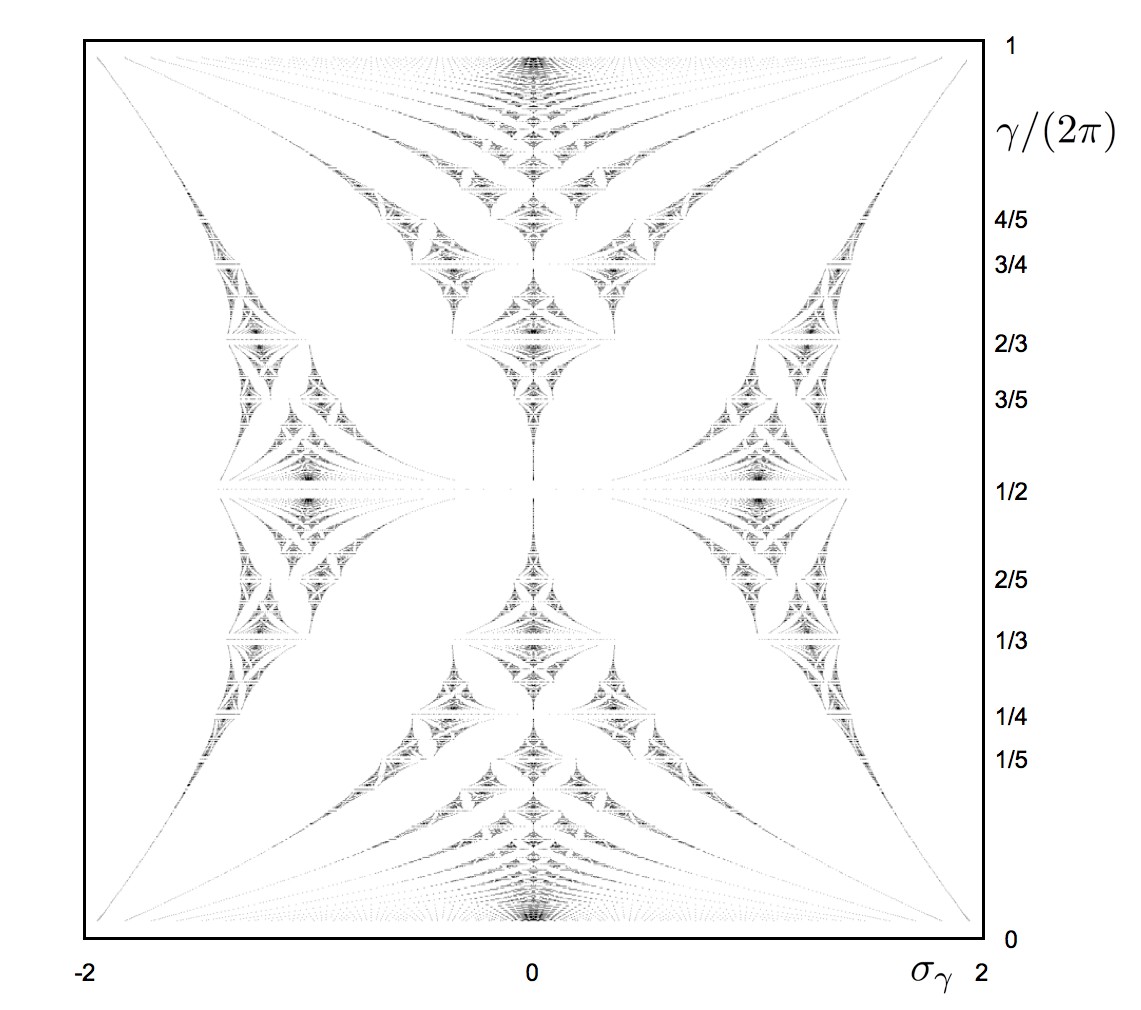} }
     	\subfloat[]{\label{square13}\includegraphics[width=.5\textwidth]{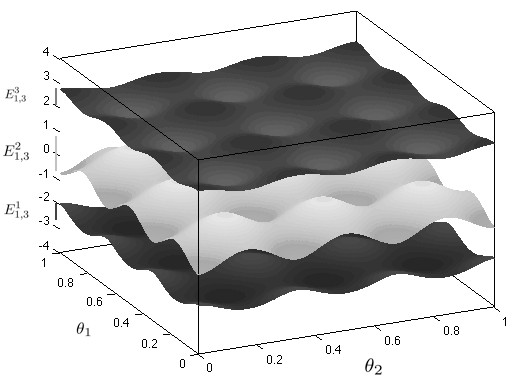} }
	\caption[The spectrum of the reduced model in the case of the square lattice]{(a) The spectrum of $P^\Box_{ \gamma}$ corresponding to the square lattice (Hofstadter's butterfly). (b) Energy bands for $\gamma= 2\pi /3$.}
 	\label{squarespectrum}
\end{figure}

Notice that $L^\Box_{ \gamma, {\theta_2+\gamma}} $ is unitary equivalent with $L^\Box_{ \gamma, {\theta_2}} $.
When $\gamma /(2\pi)$ is  irrational, the spectrum of  $L^\Box_{ \gamma,{\theta_2}}$ does not depend on ${\theta_2}$
(see \cite{HeSjHarper1}, \S 1). This is no longer the case when  $\gamma /(2\pi)$ is rational. In 1976 Hofstadter performed a
formal study of the spectrum of  $L^\Box_{ \gamma, {\theta_2}} $ as a function of $\gamma /(2\pi) \in \mathbb{Q}$
(\cite{hofstadter}). His approach suggests a fractal structure for the spectrum and leads to Hofstadter's butterfly.
The method consists in studying numerically the case $\gamma=  2\pi p/q$, with $p,q \in \N$ relative primes. Hofstadter observed
that in this case, the spectrum is formed of  $q$ bands which can only touch at their boundary. Hofstadter's butterfly
is obtained by placing in the $y$-axis of a graph the bands of  the spectrum (see  Figure \ref{spectSquare}).  Moreover,
Hofstadter derived rules for the configuration of the bands related to the expansion  of $p/q$  as continued fraction.
This configuration strongly suggests the Cantor structure  of the spectrum of $L^\Box_{\gamma,{\theta_2}} $ when
$ \gamma / (2 \pi)$ is irrational.  A longtime open problem, proposed by  Kac and Simon in the 80's and called the
``Ten Martinis problem'' (\cite{SimonXX}, Problem 4), was to prove that for irrational $  \gamma/ (2\pi)$, the spectrum of
$L^\Box_{ \gamma}$ is a Cantor set. After many efforts starting with the article of Bellissard and Simon in 1982 (\cite{BeSi}), the problem was finally solved in 2009 by Avila and Jitomirskaya (\cite{AvJi}).\\
 
In order to compute the spectrum of   $L^\Box_{ \gamma }$ for $\gamma= 2\pi  p/q$, we may use  again the Floquet theory.
Introducing the Floquet condition  $v_{n+q} =  e^{i  2\pi  \theta_1  q} v_n$, we are led to the computation of the eigenvalues
of a family (parametrized by $\theta_1$ and $\theta_2$). Denoting  $  \sigma^\Box_\gamma = \sigma(L^\Box_{ \gamma})$ we obtain  
\ben
	\sigma^\Box_\gamma  = \bigcup_{ \theta_1,\theta_2  \in [0,1]}  \sigma( M^{\Box}_{p,q,\theta_1,\theta_2}) \,,  
\een  

where
\ben
	 M^{\Box}_{p,q,\theta_1,\theta_2} =  \frac12 \left( e^{ i 2\pi   \theta_1 }K_q + e^{-i  2\pi  \theta_1 } K_q^* + e^{ i  2\pi  \theta_2 } J_{p,q} +e^{-i  2\pi  \theta_2 } J_{p,q}^* \right)  \, 
 \een

with $J_{p,q}$, $K_q$ defined in (\ref{definitpnjandk}). \\

For $1\leq k \leq q$ the $k$th band of $\sigma^\Box_\gamma $ is given by the image of
\be 
	E_{p,q} ^k:  [0,1] \times [0,1] \to \R \,,  \qquad (\theta_1,\theta_2) \mapsto  \lambda_{p,q,\theta_1,\theta_2}^k \,,
\ee

where $ \lambda_{p,q,\theta_1,\theta_2}^k$ is the $k$th eigenvalue of $ M^{\Box}_{p,q,\theta_1,\theta_2}$ 
(see  Figure \ref{square13}). \\

In  \cite{HeSjHarper1}, \S 1, it was proved that  $L^\Box_{ \gamma }$ is unitary equivalent with the pseudo-differential  operator $P^\Box_\gamma$ defined in (\ref{defGgamma}) for $p^\Box$ given in
(\ref{defpsudosqr}).  Helffer and Sj\"ostrand developed in \cite{HeSjHarper1,HeSjHarper3,HeSjHarper2} sophisticated techniques
(inspired by the work of  the physicist Wilkinson (\cite{wilkinson}) to study the operator $P^\Box_\gamma$. In particular,
they justified in various regimes the approximation for the low spectrum of $P_{h,A,V}$ by the spectrum of $P^\Box_\gamma$.   \\

When plotting $\sigma^\Box_\gamma$ as a function of $\gamma$ (see  Figure \ref{spectSquare}), we observe the following properties. 

\begin{proposition} 
\ba
	\sigma^\Box_{\gamma }  & \subset & [-2,2]   \label{squarerange} \,, \\
	\sigma^\Box_{\gamma +2\pi} &=& \sigma^\Box_\gamma  \hspace{.95cm} (\text{translation invariance}),   \label{squaretrans} \\
	\sigma^\Box_{-\gamma} &=& \sigma^\Box_{\gamma }  \hspace{.95cm} (\text{reflexion with respect to the axis $\gamma=0$}),  \label{squarereflex1}  \\
	\sigma^\Box_{2\pi-\gamma } &=& \sigma^\Box_\gamma \qquad  (\text{reflexion with respect to the axis $\gamma=\pi$}),  \label{squarereflex12}  \\
	 - e \in   \sigma^\Box_\gamma & \Leftrightarrow &   e \in  \sigma^\Box_\gamma \quad  (\text{reflexion with respect to the axis $e=0$}) \label{squarereflex2}    \,.
\ea

\end{proposition}

\textbf{Proof.} The definition of $p^\Box(x,\xi) $ together with the fact that the Weyl quantizations of $e^{ix}$, $e^{-ix}$,
 $e^{i\xi}$, $e^{-i\xi}$ are unitary operators yield (\ref{squarerange}). Property (\ref{propodd}) gives (\ref{squaretrans}).
 We obtain (\ref{squarereflex1}) noticing that 
\ben
	P^\Box_{-\gamma} = {\rm Op}^W_{-\g}\left(p^\Box(x,\xi)\right)= {\rm Op}^W_{\g}\left(p^\Box(x,-\xi)\right) 
= {\rm Op}^W_{\g}\left(p^\Box(x,\xi)\right)=P^\Box_{\gamma} \,.
\een
	
Properties (\ref{squaretrans}) and  (\ref{squarereflex1}) imply (\ref{squarereflex12}). Finally, we have that
$$p^\Box(x+\pi,\xi+\pi)=-p^\Box(x,\xi)$$
so $P^{\Box}_{\g}$ and $-P^{\Box}_{\g}$ are conjugate by the unitary
operator $u\mapsto e^{\frac{i}{\gamma}\pi\cdot}u(\cdot-\pi)$. This yields (\ref{squarereflex2}).\\ \qed

 \subsection{The triangular lattice}  

\begin{figure}  [H]
  \centering
  	\subfloat[]{\label{spectTriangular}\includegraphics[width=.425\textwidth]{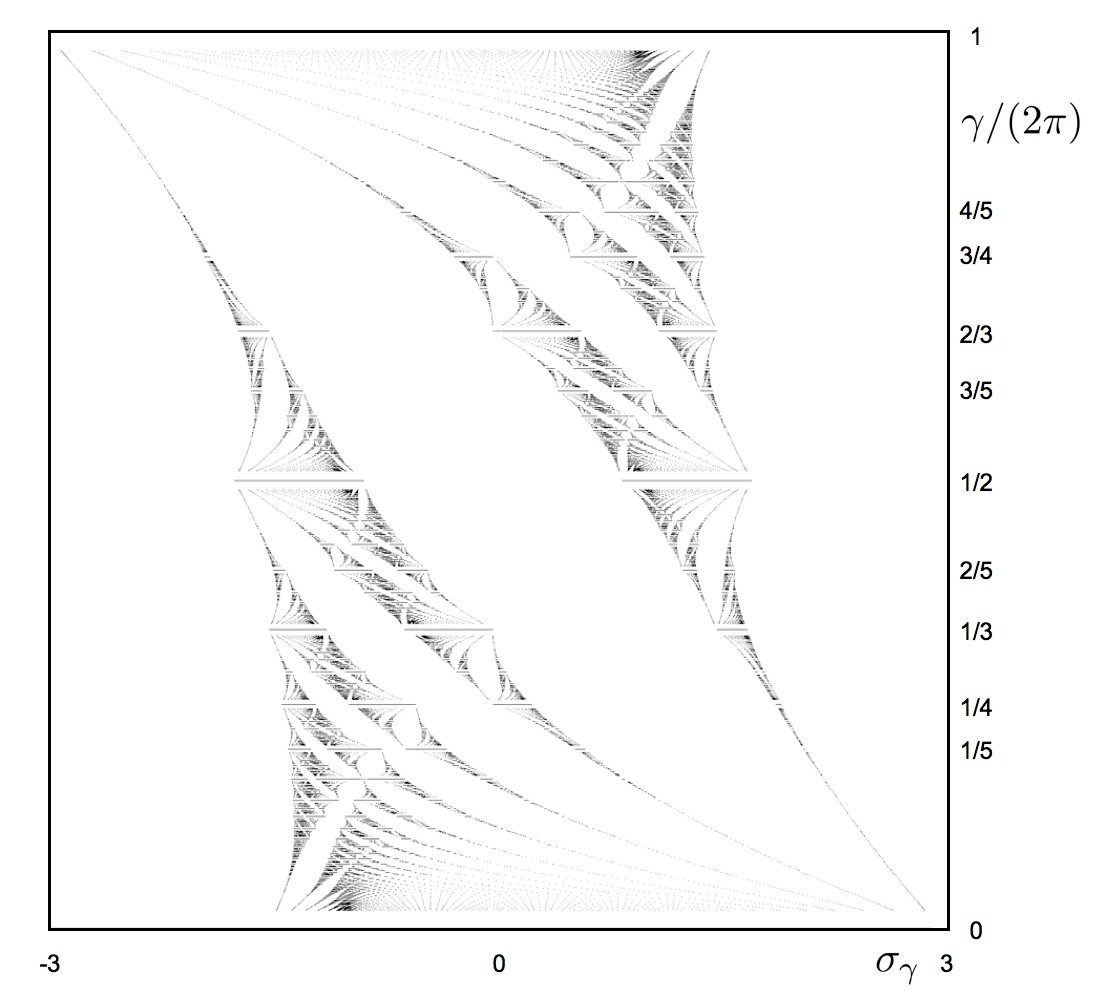} } 
     	\subfloat[]{\label{spectHexagonal}\includegraphics[width=.425\textwidth]{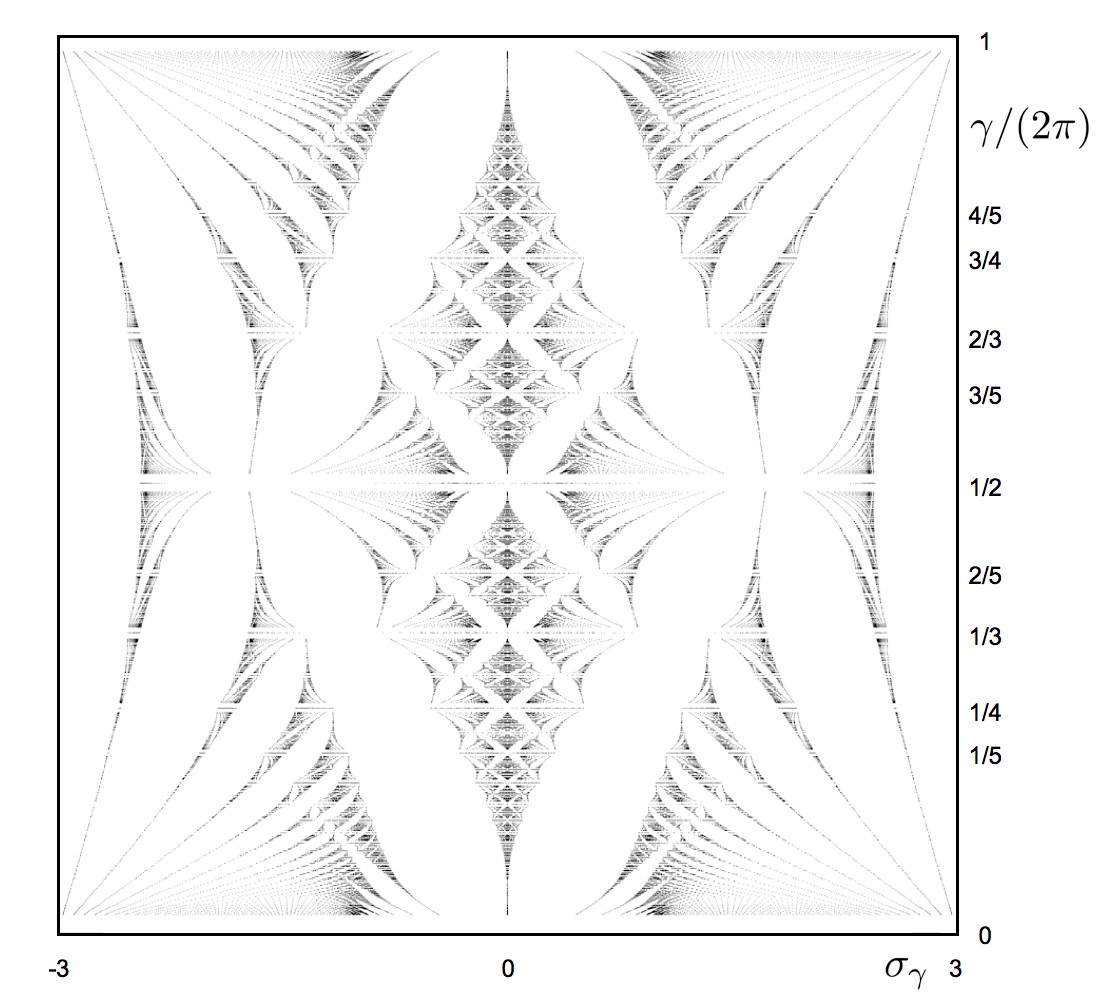} }
	\caption[Hofstadter's butterfly for the triangular and the hexagonal lattices]{Hofstadter's butterfly for (a) the triangular and (b) the hexagonal lattices.}
 	\label{trihexspectrum}
\end{figure}

The triangular lattice\footnote{We note that the triangular and hexagonal lattices are sometimes respectively called hexagonal
and honeycomb lattices.} is the Bravais lattice associated with the basis $\{(1,0), (\nicefrac 12, -\nicefrac {\sqrt{3}}2) \}$.
Each point of the lattice has 6 nearest neighbors for the Euclidean distance. This case was studied by Claro and Wannier in
\cite{ClaroWannier}. These authors exhibit an analogous structure to the case of the square lattice.
In the case $\gamma= 2\pi  p/q$, with $p,q \in \N$ relative primes, the spectrum is formed of  $q$ bands which can only
touch at their boundary  (see Figure \ref{spectTriangular}). In \cite{kerdelhue1}, the first author studied rigorously the operator
$P_{h,A,V}$ in this case. He justified the reduction to  the pseudo-differential  operator $P_\gamma^\triangle$ defined in
(\ref{defGgamma}) with $p^\triangle$ given in (\ref{triangularsymbol}). \\

As in the case of the square lattice discussed before, when $\gamma=2\pi p/q$ the spectrum can be computed by considering the 
family of matrices in $M_{ q}(\C)$ defined by
\ban
	M^{\triangle} _{p,q,\theta_1,\theta_2} & =& \frac12 \left( \ e^{i   2\pi \theta_1} K_q + e^{-i   2\pi \theta_1} K_q^* + e^{i   2\pi \theta_2}J_{p,q} + e^{-i   2\pi \theta_2} J_{p,q}^* \right. \\
	&&  \left. +  e^{- i  \pi \nicefrac pq}  e^{i   2\pi ( \theta_1+ \theta_2)} J_{p,q} K_q  + e^{- i  \pi \nicefrac pq}  e^{ -i 2\pi( \theta_1+\theta_2 )}  J_{p,q}^* K_q^*  \right)\,.
\ean
 
Let $\sigma^\triangle_\gamma$ be the spectrum of $P^\triangle_\gamma$.  When plotting $\sigma^\triangle_\gamma$  as a
function of $\gamma  $ (see  Figure \ref{spectTriangular}), we observe the following properties. 

\begin{proposition}  
\ba
	 \sigma^\triangle_{\gamma }  & \subset & [-3,3] \,, \label{trirange}  \\
	\sigma^\triangle_{\gamma +4\pi} &=& \sigma^\triangle_\gamma  \hspace{.95cm} (\text{translation invariance}),  \label{tritrans} \\
	\sigma^\triangle_{-\gamma} &=& \sigma^\triangle_{\gamma }  \hspace{.95cm} (\text{reflexion with respect to the axis $\gamma=0$}), \label{trireflex1}  \\
	  e \in   \sigma^\triangle_{2\pi-\gamma} & \Leftrightarrow &   -e \in  \sigma^\triangle_{\gamma } \qquad  (\text{reflexion with respect to the point $(0,\pi)$}) \label{triantitrans}   \,. 
\ea

\end{proposition}

\textbf{Proof.} The definition of $p^\triangle(x,\xi) $ together with the fact that the Weyl quantizations of $e^{ix}$, $e^{-ix}$,
$e^{i\xi}$, $e^{-i\xi}$, $e^{i(x-\xi)}$ and $e^{-i(x-\xi)}$ are unitary operators yield (\ref{trirange}). 
 Property (\ref{tritrans}) comes from Proposition \ref{propositionsym}.  We have that
\ben 
	 {\rm Op}^W_{-\g}\left(p^\triangle(x,\xi)\right)= {\rm Op}^W_{\g}\left(p^\triangle(x,-\xi)\right)
= {\rm Op}^W_{\g}\left(\cos x+\cos \xi+\cos(x+\xi)\right) \,.
\een

The application $(x,\xi)\mapsto (x,-x+\xi)$ is linear symplectic so $ {\rm Op}^W_{\g}\left(\cos x+\cos \xi+\cos(x+\xi)\right)$ 
and 
\ben
	 {\rm Op}^W_{\g}\left(\cos x+\cos (-x+\xi)+\cos(x+(-x+\xi))\right)= {\rm Op}^W_{\g}\left(p^\triangle(x,\xi)\right)
\een	 

are unitary equivalent. This yields (\ref{trireflex1}). Proposition  \ref{propositionsym} implies that
$ e \in   \sigma^\triangle_{2\pi-\gamma} $ if and only if  $-e \in  \sigma^\triangle_{-\gamma } $. This, together with (\ref{trireflex1}), yield (\ref{triantitrans}). \\ \qed

 \subsection{The hexagonal lattice} 

The hexagonal lattice is not a Bravais lattice, but is a discrete subset of $\mathbb R^2$ invariant under
the rotation of angle $\pi/3$ and translation along a triangular lattice, and containing two points per fundamental domain of this lattice. 
Each point of the lattice has 3 nearest neighbors. This case was also  rigorously  studied by the first author
in \cite{kerdelhue1} and \cite{kerdelhue2}. We remark that this configuration corresponds to a charged particle in a graphene
sheet submitted to a  transverse magnetic  field (\cite{polyMontambaux}, \S 6). This case acquired a new interest after the
2010 Nobel Prize in Physics awarded to Geim and Novoselov for their experiments involving graphene
(\cite{Geim,nobel2010,Novoselov}). In the case of a hexagonal lattice, thefirst author justified the reduction to a
pseudo-differential operator $P_\gamma^\hexagon$ defined in (\ref{defGgamma}) with $p^\hexagon$ given in (\ref{hexagonalsymbol}). \\

In the case  when $\gamma=2\pi p/q$, the spectrum can be numerically computed by diagonalizing  the hermitian matrices 
in $M_{2q}(\C)$ defined by
\ben %
M^{\hexagon}_{p,q,\theta_1,\theta_2} = \left(\begin{array}{c|c} 
	&       \\
	  0_q   &  I_q +  e^{i  {\theta_1}} K_q +  e^{-i  {\theta_2}} J^*_{p,q}   \\
	  &       \\
	\hline 
	  &       \\
	   I_q +  e^{-i  {\theta_1}} K_q^* +  e^{ i  {\theta_2}}J_{p,q}   &   0_q       \\
	  &    
	\end{array}\right) \,.
\een

Let  $\sigma^\hexagon_{\gamma}$ be the spectrum of  $P_\gamma^\hexagon$.  When plotting $\sigma^\hexagon_{\gamma}$ as a function of $\gamma  $ (see  Figure \ref{spectHexagonal}), we observe the following properties.

\begin{proposition}  \label{symhex}
\ba
	 \sigma^\hexagon_{\gamma }  & \subset & [-3,3] \,, \label{hexrange} \\
	\sigma^\hexagon_{\gamma +2\pi} &=& \sigma^\hexagon_{\gamma}  \hspace{.95cm} \text{ (translation invariance), }  \label{hextrans} \\
	\sigma^\hexagon_{-\gamma} &=& \sigma^\hexagon_{\gamma}  \hspace{.95cm} \text{ (reflexion with respect to the axis $\gamma=0$}), \label{hexreflex1}  \\
	\sigma^\hexagon_{2\pi-\gamma } &=& \sigma^\hexagon_{\gamma} \qquad  \text{ (reflexion with respect to the axis $\gamma=\pi$}),  \label{hexreflex12}  \\
	 - e \in   \sigma^\hexagon_\gamma & \Leftrightarrow &   e \in  \sigma^\hexagon_{\gamma}
\quad  (\text{reflexion with respect to the axis $e=0$}) \label{hexreflex2}
	 \,.
\ea

\end{proposition}

\textbf{Proof.}  We obtain (\ref{hexrange}) observing that 
\ben 
p^{\hexagon}(x,\xi)=\left(\begin{array}{cc}
                            0&1\\
1&0
                           \end{array}\right)+
\left(\begin{array}{cc}
                            0&e^{ix}\\
e^{-ix}&0
                           \end{array}\right)+
\left(\begin{array}{cc}
                         0&   e^{i\xi}\\
e^{-i\xi}&0
                           \end{array}\right)
\een

and that  the Weyl quantizations of 
\ben
	\left(\begin{array}{cc}
                            0&1\\
1&0
                           \end{array}\right) \,, \qquad \left(\begin{array}{cc}
                            0&e^{ix}\\
e^{-ix}&0
                           \end{array}\right) \quad \text{ and } \quad  
 \left(\begin{array}{cc}
                           0 & e^{i\xi}\\
e^{-i\xi}&0
                           \end{array}\right) 
\een                           
are unitary operators. Property (\ref{hextrans}) comes from Proposition \ref{propositionsym}. Notice that 
\ben 
 {\rm Op}^W_{-\g}\left(p^{hex}(x,\xi)\right)= {\rm Op}^W_{\g}\left(p^{hex}(x,-\xi)\right)
= {\rm Op}^W_{\g}\left(\begin{array}{cc}
                            0&1+e^{ix}+e^{-i\xi}\\
1+e^{-ix}+e^{i\xi}&0
                           \end{array}\right)  \een
                           
and let $\Gamma$ be the operator defined by $\Gamma u(x)=\overline{u(-x)}$. It is classical and easy to check that if $q$ is a symbol,
$\Gamma\, {\rm Op}^W_{\g}q(x,\xi)\,\Gamma= {\rm Op}^W_{\g}\overline{q(-x,\xi)}$. This gives  
\ben
	\Gamma\, {\rm Op}^W_{\g}\left(\begin{array}{cc}
                            0&1+e^{ix}+e^{-i\xi}\\
1+e^{-ix}+e^{i\xi}&0
                           \end{array}\right)\,\Gamma
= {\rm Op}^W_{\g} p^{hex} (x,\xi) \,,
\een

which yields (\ref{hexreflex1}). Property (\ref{hexreflex12}) follows from (\ref{hextrans}) and (\ref{hexreflex1}).
Finally, noting that 
 \ben
 \left(\begin{array}{cc}
            -1 & 0\\
0 & 1
           \end{array}\right)
p^{hex}(x,\xi)
\left(\begin{array}{cc}
            -1 & 0\\
0 & 1
           \end{array}\right)
=-p^{hex}(x,\xi)  \, 
\een

we obtain
 \ben
 \left(\begin{array}{cc}
            -1 & 0\\
0 & 1
           \end{array}\right)
P^{hex}_{\g}
\left(\begin{array}{cc}
            -1 & 0\\
0 & 1
           \end{array}\right)
=-P^{hex}_{\g} \,,
\een 

which yields (\ref{hexreflex2}).\\  \qed

 \section[The kagome lattice]{The kagome lattice} \label{lattice} 

\subsection{The group of symmetries of $\Gamma$}
\label{deflabel}

We now study the properties of the kagome lattice and its group of symmetries $\mathcal{G}$.\\

For $\a \in \Z^2$ we set $t^\a = t_1^{\a_1} t_2^{\a_2} $. We then have  
\be \label{rtaualpha}
	rt^{\alpha}   = t^{\kappa(\alpha)} r \,,
\ee 

where $\kappa $ is given in (\ref{defkappa}). We also notice that  
\be \label{relationsskappa}
	\kappa^6=id_{\Z^2}  \quad \text{ and } \quad \kappa(\a) \wedge \kappa(\b) = \a \wedge \b \,,
 \ee

where $\wedge$ is the cross product  $\a \wedge \b =  \a_1\b_2 -  \a_2\b_1 $. Then we easily obtain 

\begin{proposition} \label{propGG}
The kagome lattice is invariant by the maps in $\mathcal{G}$ and for every $m,n \in \Gamma$  there exists $g \in \mathcal{G}$ such that $g(m)=n$.
\end{proposition}

\subsection{Construction of kagome potentials}   
\label{construction}

We call $V : \R^2 \to \R$ a \textit{kagome potential} if it satisfies Hypothesis \ref{hypothesisPotV}.  It is rather easy to define such a potential, but more interesting is to give explicit examples in the class of trigonometric polynomials, which leaves open the possibility  to realize experimentally these potentials with lasers (see for example \cite{damski} and \cite{santos}). \\  

%

\begin{figure} [H]
	\center  \includegraphics[width=.6\textwidth]{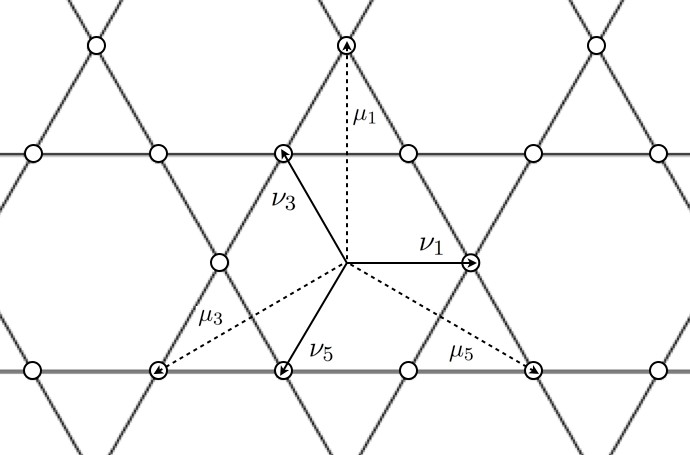}  
	\caption[Vectors $\nu_1$, $\nu_3$, $\nu_5$,  $\mu_1$, $\mu_3$ and  $\mu_5$]{Vectors $\nu_1$, $\nu_3$, $\nu_5$,  $\mu_1$, $\mu_3$ and  $\mu_5$.}
\label{mus}
\end{figure}

Remembering the definitions of the vectors $\nu_j$ from (\ref{defnu}), we denote by $\nu^\perp$ the vector deduced from 
$\nu$ by a rotation of $\nicefrac \pi2$ and for $j \in\{1,3,5 \}$ we define (see Figure \ref{mus})
\be \label{fedmuj}
	\mu_j = \sqrt{3} \,   \nu_j^\perp \,.
\ee

%
%
%
%

For $j =1,3,5$ we set $\phi_j=3\pi/2$ and define the potentials $V_j : \R^2 \to \R$ as
\be \label{defVjkagpot}
	V_j(x) = \left[ \cos \left(  x \cdot \pi \mu_j + \phi_j \right) + 
2  \cos \left(  \frac{x \cdot \pi  \mu_j + \phi_j}3 \right) \right]^2  \,,
 \ee 

and $\tilde V$ as
\be \tilde V=V_1 + V_3 +V_5 \,.\ee

A straightforward computation gives 

\begin{proposition}
The function 
\be \label{defVkagpot} 
	V =-\tilde V+ \|\tilde V\|_{\infty} \,,
\ee 
 satisfies (\ref{symmV}) and (\ref{nodegV}) and has local minima at the points of the kagome lattice.
\end{proposition}

\begin{figure}[H] 
  \centering
  	\subfloat[]{\label{superlatticep2}\includegraphics[width= .3\textwidth]{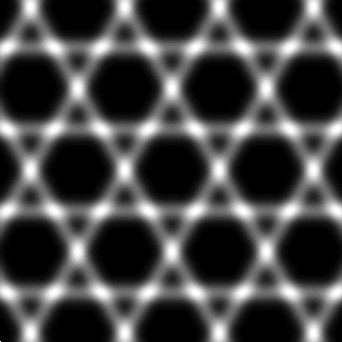}}     
	\hspace{.25cm}    
  	\subfloat[]{\label{superlatticep4}\includegraphics[width=.3\textwidth]{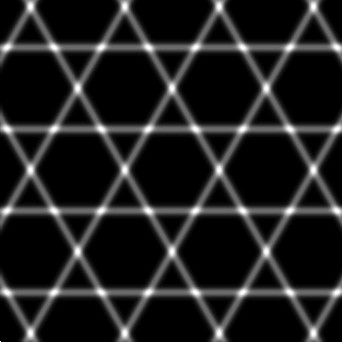}} 
	\hspace{.25cm}
	\subfloat[]{\label{superlatticep100}\includegraphics[width=.3\textwidth]{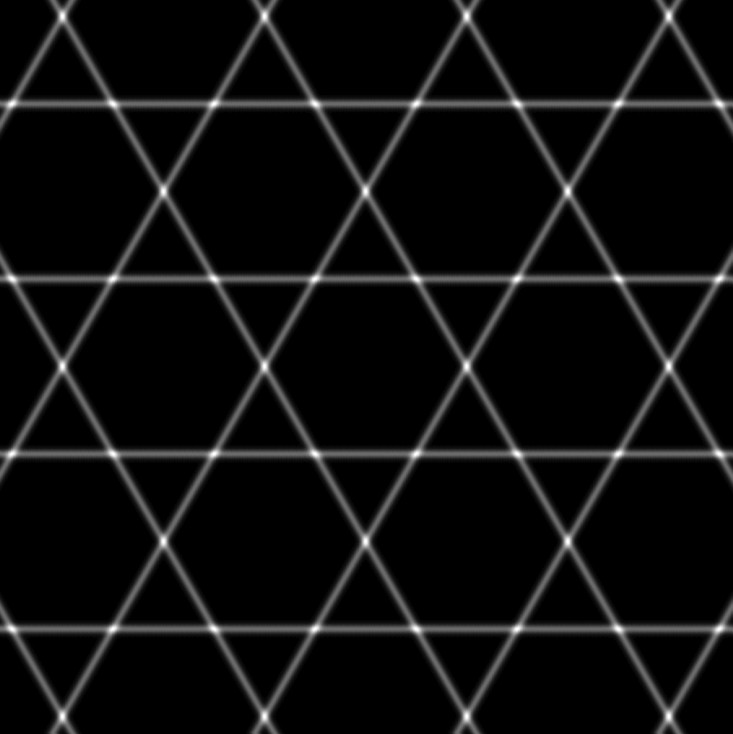}}
	 \caption[kagome potential for $p \in \{2,10,30,50 \}$ ]{The Kagome potentials  $V =\|V_1^{\nicefrac p2}+V_3^{\nicefrac p2}+V_5^{\nicefrac p2}\|_\infty-(V_1^{\nicefrac p2}+V_3^{\nicefrac p2}+V_5^{\nicefrac p2})$ for (a) $p=2$, (b) $p=10$, (c) $p=40$ with $V_j $ given by (\ref{defVjkagpot}) for $j \in \{1,3,5\}$. Smaller values are represented by darker colors.}
   \label{superlattice}
\end{figure}

\begin{remark}
	Our numerical computations (see Figure \ref{superlatticep2}) show that the condition (\ref{minimaV}) is verified but we do not have a mathematical proof.
\end{remark}

\begin{remark}
	We notice that the potential defined by (\ref{defVkagpot}) with 	
\ben
	V_j(x) = \left[ \cos \left(  x \cdot \pi \mu_j + \phi_j \right) + 2  \cos \left(  \frac{x \cdot \pi  \mu_j 
+ \phi_j}3 \right) \right]^p \,, \quad p \in 2\mathbb N \,,
 \een 
 is also a kagome potential (see Figure \ref{superlattice}).  When $p$ goes to $+\infty$, we observe that the minima are very well localized at the points of  $\Gamma$. This could be an advantage for verifying theoretical  assumptions for an accurate semi-classical analysis of the tunneling effect between wells in the next section, but $p$ large is not experimentally reasonable.    
\end{remark}

\begin{remark}
	Considering any Bravais lattice with three points by periodicity cell, we are led to the same situation, but the kagome lattice have a much richer structure.  
\end{remark}

\section{The Schr\"odinger magnetic operator on $L^2(\R^2)$} \label{operator}   

\subsection{The Schr\"odinger magnetic operator}

 We recall that we start from $P^0_{h,A,V}$ defined in (\ref{ScrhoPot3}). Since we have assumed  $V\geq0$, the operator is
semi-bounded on $C_0^\infty(\R^2) $ and there is an unique selfadjoint extension in $L^2(\R^2)$, which can be obtained as the 
Friedrichs extension of $P^0_{h,A,V}$ (see for example \cite{HeBook}). It can be proved that the domain of $P_{h,A,V}$ is 
given by
\be
	\mathcal D(P_{h,A,V}) = \Big\{  u \in L^2(\R^2) \,;\, P_{h,A,V}\,u \in L^2(\R^2)   \Big \} \,.
\ee

\subsection{Quantization of $\mathcal{G}$}

The use of the symmetries in the case of the square, triangular and hexagonal lattices was crucial in \cite{HeSjHarper1} and \cite{kerdelhue1}. In order to take advantage of the properties of the kagome lattice, we need to quantify the elements of $\mathcal{G} $, that is, to associate which each element of $\mathcal{G} $ an unitary transformation in $L^2(\R^2)$, which respects the domain and commutes with $P_{h,A,V}$. These operators will be used later to study  the low lying spectrum of $P_{h,A,V}$. We note that the quantization of the translations $T_j$ was introduced by Zak in \cite{zak}. We also mention the work of Helffer and Sj\"ostrand (\cite{HeSjLecture}, pages 147-148) who studied the case of constant magnetic field in arbitrarily dimension (see also Bellissard (\cite{Be1}), Cartier (\cite{Cartier}) and Zak.   \\

Since the symmetries of the kagome lattice are dictated by those of the triangular lattice, we will use the construction of
the first author in Section 1 of \cite{kerdelhue1}. We explain in the following the main ideas.

\subsubsection{Quantization of the rotation and the translations}

We now quantify the rotation $r$ and the translations $t_j$. We notice that for every $g \in \mathcal{G} $ the $1$-form $A-gA$ 
is closed and in fact it is exact. Indeed, by assumption (\ref{symmsB}),  
\be \label{exactform}
	d(A-gA)=dA-g \,dA= B - g B = 0 \,.
\ee

Hence, there is a real smooth function $\phi_g$, defined up to a constant, such that 
\be
\label{orange}
	A-gA=d\phi_g \,.
\ee

Later, we will use this freedom of choice of the constants to obtain simple commutation properties. \\

 We may then quantize  $g \in \mathcal{G} $ by the operator $T_g$, defined on $C^\infty_0(\R^2)$ by
\be \label{defquant1}
	(T_gu)(x)=e^{\frac ih \phi_g(x)} \, u(g^{-1}(x)) \,,
\ee

 where $\phi_g$ is the a real function associated with $g$ by (\ref{orange}).

\begin{lemma} \label{dosh}
	For any $g \in \mathcal{G}$, the operator $T_g$ is unitary on $L^2(\R^2)$ and commutes with $P_{h,A  ,V}$. 
\end{lemma}

\textbf{Proof.} For the first assertion a simple computation gives 
\be
\label{worm}
		T_g^{-1} = e^{-\frac{i}h  g^{-1} \phi_g }\, g^{-1} = T_g^* \,.
\ee

We have the equality between 1-forms  
\ban
	(-ihd-A) T_g u &=& e^{i \nicefrac{\phi_g}h}  \, (  (d\phi_g) \, g u -i h \, g  (d  u)  - A\, (g  u) )  \,,
\ean

so using (\ref{orange}) we get 
\ban
	(-ihd-A) T_g u  &=&  e^{i \nicefrac{\phi_g}h}   \, (-ih \, g  (d  u)  - (g A) \,( g  u ) ) \\
				     &=&  e^{i \nicefrac{\phi_g}h}  \, g ( -ih d u  - A u  ) \\
				     &=& T_g (-ihd-A) u \,,
 \ean

which gives the lemma. \\ \qed 

\subsubsection{Definition of the magnetic rotation and translations} 

For $j=1,2,3 $ we define the magnetic translations
\be
	T_j = e^{\frac ih \phi_j} \, t_j \,,
\ee

where $\phi_j$ is the real function associated with $t_j$ by (\ref{orange}) with $g=t_j $.  \\

The inverse of $T_j $ is given by (\ref{worm}) and is also a magnetic translation. For  $j=4,5,6 $ we then define 
\be \label{defTransinv}
	T_j = T_{j+3}^{-1} \,.
\ee 

We also define the magnetic rotation
\be
	F =e^{\frac ih f } \,  r^{-1} \,,
\ee

where $f$ is the real function associated with $g=r^{-1}$ by (\ref{orange}).

\begin{remark}
	We need a convenient choice of $T_{g_1g_2}$ in order to be able to compare with $T_{g_1} \circ T_{g_2}$ 
(see (\ref{definitionTalpha}) below). Hence, we will only use the previous construction for $r $ and $t_j$, $j \in \Zm$\,. 
\end{remark}

\begin{remark} \label{remark1}
 In the case of a constant magnetic field, choosing the gauge $A(x_1,x_2)=\frac B2 (-x_2,x_1)$ we have 	
\ben
	f(x) =  f_0 \qquad \text{ and } \qquad  \phi_j(x) = -  \frac B2 \, x \wedge(2 \nu_j)   +c_j \,,  
\een 
 where $f_0$ and  $c_j$, $j=1,2,3$, are arbitrary constants. 
\end{remark}

\subsubsection{Commutation rules}

We now show how a good choice of the constants appearing in the definition of $f$ and $\phi_j$ lead to nice commutation rules for the operators $F$ and $T_j$. 

\begin{proposition} \label{magneticflux}

\textbf{(i)} The flux of $B$  through a fundamental domain $\mathcal V$ of $\Gamma_\triangle$ does not depend on the basis 
chosen. We write
\ben
	\gamma = \frac1h \int_{ \mathcal V } d\omega_A \,.
\een

\textbf{(ii)} We have
\be \label{flux}
	T_j T_{j+1} = e^{ i \gamma}T_{j+1}T_j \,.
\ee  

\textbf{(iii)} There are unique $\phi_1$, $\phi_2$, $\phi_3$ and $f$ such that
\ba 
	F^6 &=& Id_{L^2(\R^2)} \,, \label{relaF6} \\
	T_jF &=& FT_{j+1}  \,, \label{relaFT}
\ea  

and for these $\phi_1$, $\phi_2$ and $\phi_3$ we have
\be \label{TjTj2}
	T_j T_{j+2} = e^{i \frac{  \gamma}{2 }}T_{j+1}  \,.
\ee 

 \end{proposition}

\textbf{From now on, we choose $\phi_1$, $\phi_2$, $\phi_3$ and $f$ in the definition of $T_1$, $T_2$, $T_3$ and $F$ such that (\ref{relaF6})-(\ref{TjTj2}) are satisfied.} \\

\begin{remark}  \label{remark2}
	In the case of a constant magnetic field, choosing $A(x_1,x_2)=\frac B2 (-x_2,x_1)$ we verify  	
\ben
	f_0   = c_1 = c_2 = c_3 =0 \,.
\een	
 	 
\end{remark}

\textbf{Proof of Proposition \ref{magneticflux}.} The translations $t_j$ commute between them so we have   
\ben
	T_j T_{j+1} = e^{  \frac ih  \{  \phi_j+t_j \phi_{j+1}-\phi_{j+1}-t_{j+1} \phi_j   \}} T_{j+1} T_j \,.
\een

 After (\ref{exactform}), the expression between the brackets here before is a constant that we note $\eta_j$. Using 
(\ref{orange}) we compute
\ban
	(\phi_j - t_{j+1}  \phi_j )(x) &=& \int_{[x-2\nu_{j+1},x]} d\phi_j\\
	 &=& \int_{[x-2\nu_{j+1},x]} ( A-t_j  A ) \\
	&=&   \int_{[x-2\nu_{j+1},x]}  A +   \int_{[x-2\nu_{j},x-2\nu_{j}-2\nu_{j+1}]}    A \,,
 \ean  

where $[x,y]$ denotes the path $[0,1] \ni s  \mapsto   (1-s)x+sy$. Similarly,
\ban
	(t_j \phi_{j+1} - \phi_{j+1} )(x) &=& \int_{[x,x-2\nu_j]} A +  \int_{[x-2\nu_j- 2\nu_{j+1},x- 2\nu_{j+1}]}   A \,.
\ean  

Hence, Stokes theorem yields 
\be \label{stokes}
	\eta_j  = \int_{\mathcal{V}_{j,j+1}} d\omega_A  \,,
\ee 

where $\mathcal{V}_{j,j+1}$ is a cell of periodicity of the lattice generated by $2\nu_j$ and $2\nu_{j+1}$ with vertex $x$, $x-2\nu_{j+1}$, $x-2\nu_{j}$ and $x-2\nu_{j}-2\nu_{j+1}$. After (\ref{symmsB}) the magnetic field $B dx_1 \wedge dx_2$ is invariant by $r$, so the value of the $\eta_j$ do not depend on $j \in \Zm$ and we have $\gamma=\eta_1/h$. We have proved \textbf{(i)} and \textbf{(ii)}. \\

Since $r^6=id_{\Z^2}$ we have
\be \label{F3}
	F^6= e^{ \frac ih \{  f+rf+\cdots+ r^5f \} }  \,.
\ee 

After (\ref{exactform}) the expression between the brackets here before is a constant. Hence, choosing an appropriate constant in the definition of $f$ we obtain (\ref{relaF6}). \\

Using (\ref{rtaualpha}) and (\ref{defTransinv}) we have
\ban
	T_1F &=& e^{\frac ih \{ \phi_1 + t_1  f -  f - r^{-1} \phi_2 \}} F  T_{2} \,, \\
	T_2F &=& e^{\frac ih \{ \phi_2 + t_2  f -  f - r^{-1} \phi_3 \}} F  T_{3} \,, \\
	T_3F &=& e^{\frac ih \{ \phi_3 + t_3  f -  f + t_3 r^{-1} \phi_1 \}} F  T_4 \,.
\ean 

As before, using (\ref{exactform}) the expressions between the brackets in the previous equalities are constants.
If we respectively add $a_1$, $a_2$ and $a_3$ to $\phi_1$, $\phi_2$ and $\phi_3$, the expressions between the brackets are
respectively modified by $ a_1- a_2$, $ a_2-a_3$ and $ a_1+a_3$. Hence, there exist $a_1$, $a_2$ and $a_3$ such that (\ref{relaFT}) is satisfied for $j=1,2,3$. Since $T_{j+3} = T_j^{-1}$, (\ref{relaFT}) also holds for  $j=4,5,6$. \\

Again, we have that $\{ \phi_j+t_j \phi_{j+2}-\phi_{j+1} \}$ is a constant that we call $c$, so 
\be \label{F33}
	T_j T_{j+2}= e^{i\frac ch} \, T_{j+1} \,.
\ee 

 Using the conjugation by $F$ and (\ref{relaFT}), we then obtain 
$  e^{i\frac ch} \, T_{j+2} = T_{j+1} T_{j}^{-1} $, which gives  
\be
	T_j T_{j+2}= 	e^{i\frac{2c}h} T_{j+2} T_{j} \,.
\ee

 The proof of (\ref{flux}) also applies when taking $T_j$ and $T_{j+2}$ instead of $T_j$ and $T_{j+1}$, so we have $T_j T_{j+2} = e^{i\frac {\eta_1} h} T_{j+2} T_{j}$. Thus, $  2c/h \equiv  \eta_1/ h  \,\, [ 2\pi]$, which gives $ (c -   \eta_1/2)/h \in   \pi  \Z$.  Since $c$ and $\eta_1$ do not depend on $h$, we derive that necessarily  $c =   \eta_1/2 $, which yields   (\ref{TjTj2}). We have proved \textbf{(iii)}.   \\ \qed \\

Now, for $\a \in \mathbb{Z}^2$ we define the magnetic translations 
\be \label{definitionTalpha} 
	T^\a = e^{-i\frac{\gamma}{2 } \a_1 \a_2 } \, T_1^{\a_1}T_2^{\a_2} \,.
\ee

We obtain the following relations (see also \cite{kerdelhue1}, pages 15-16): 

\begin{proposition}
\label{tecito} 
For every $\a,\b \in \mathbb{Z}^2$, 
\ba 
	(T^{\a})^{-1} &=& T^{-\a} _,,\\
	T^{\a} T^{\beta} &=& e^{i\frac{\gamma  }{2 } \a \wedge \beta} \, T^{\a+\beta} \label{Talphaplusbeta} \,,\\	
	FT^{\a}&=& T^{\kappa^{-1}(\a)}F\label{TalphaF}  \,.
\ea		
\end{proposition}
 
\textbf{Proof of Proposition \ref{tecito}.} Using (\ref{flux}) we have 
$T_1^{\a_1}T_2^{\a_2}=e^{i \a_1\a_2 \gamma}T_2^{\a_2}T_1^{\a_1}$, which gives
\ban
	(T^{\a})^{-1}  &=&   e^{ \frac  {i \gamma }2   \a_1 \a_2}  \,\,T_2^{-\a_2}T_1^{-\a_1} \\
				   &=& e^{ \frac  {i \gamma }2  (\a_1 \a_2 - 2 \a_1 \a_2)}  \,\, T_1^{-\a_1} T_2^{-\a_2}\\
				   &=& T^{-\a} \, 
\ean  

and
\ban
	T^{\a}T^{\beta}  &=&   e^{-  \frac  {i \gamma }2 ( \a_1 \a_2+ \beta_1\beta_2)} \,\,T_1^{\a_1}T_2^{\a_2} T_1^{\beta_1}T_2^{\beta_2} \\
				   &=& e^{-  \frac  {i \gamma }2  ( \a_1 \a_2+ \beta_1\beta_2 + 2 \a_2 \beta_1)} \,\,  T_1^{ \a_1 + \b_1} T_2^{\a_2+\b_2} \\
				   &=& e^{  \frac  {i \gamma }2  ( \a \wedge \beta)}\,\,   T^{\a + \b} \,.
\ean 

Using (\ref{TjTj2}) we have $T_1^{ \a_1}T_3^{\a_1} =e^{ \frac  {i \gamma }2 \a_1^2 }T_2^{ \a_1}$. Hence, after 
(\ref{defTransinv}) and (\ref{relaFT}) we get 
\ban
	F T^{\a}   &=&   e^{- \frac  {i \gamma }2  \a_1 \a_2} \,\, F T_1^{\a_1}T_2^{\a_2}  \\
		        &=& e^{- \frac  {i \gamma }2  \a_1 \a_2} \,\,  T_3^{-\a_1}  T_1^{\a_2} F \\
		        &=& e^{- \frac  {i \gamma }2  ( \a_1 \a_2 + \a_1^2)} \,\,  T_2^{-\a_1} T_1^{\a_1+ \a_2} F \\
		        &=& e^{- \frac  {i \gamma }2  (\a_1 + \a_2  ) (-\a_1)} \,\,  T_1^{\a_1+\a_2} T_2^{-\a_1} F \\
		        &=& T^{k^{-1}(\a)}F \,.
\ean 

 \qed

\subsection{The harmonic approximation}  \label{harappsection}

Here we recall a result from \cite{HeSjPise} about the semiclassical analysis of the bottom of the spectrum of a Schr\"odinger 
operator with  magnetic field in the case when the electric potential $V$ has a unique non degenerate well at a point $m$.\\

The theory of the harmonic approximation was initially introduced for a Schr\"odinger operator without magnetic field in 
\cite{HeSjWsm1} and \cite{SiIHP} and can be extended to the magnetic case. More precisely,  the harmonic approximation consists
 in replacing the potential $V$ by its quadratic approximation at $m$ and the magnetic field  by its value at $m$, that is the
 magnetic potential by its linear part at $m$. This reads:  
\be \label{harmagapprox}
	P_{har,h,B} = h^2D_{x_1}^2 + (hD_{x_2} - B  x_1)^2 + \frac12 \langle {\rm Hess}\,V(m) \,x,x\rangle \, 
\ee 

with $B=B(m)$. The following result is classical and can for example be found in \cite{HelfferSyrie}: 
 
 \begin{proposition} \label{harapp0} 
Assume that ${\rm Hess}\,V(m)>0$. The spectrum of the operator $P_{har,h,B}$ defined in (\ref{harmagapprox}) is discrete. The first 
eigenvalue is simple and given by 
\ben
	\lambda_{ har, h,B } = h \sqrt{ \lambda^2_{1,0}+ B^2 } \,,
\een

where $\lambda_{ 1,0}= (  \sqrt{\lambda_1}  + \sqrt{ \lambda_2}  )/\sqrt{2} $ is the first eigenvalue of $P_{har,1,0}$ and 
$ \lambda_1$,$ \lambda_2$ are the eigenvalues of  ${\rm Hess}\,V(m)$.

 \end{proposition}
 
\textbf{Proof.} Possibly after changes of coordinates and gauge,  $P_{har,h,B}$ is written
\ben 
	P_{har,h,B} = h^2D_{x_1}^2 + (hD_{x_2} - B  x_1)^2 + \frac{ \lambda_1}{2} x_1^2 + \frac{ \lambda_2}{2} x_2^2 \,.
 \een

A partial Fourier transform in the second variable leads to the operator 
\ben h^2D_{x_1}^2 + (h\xi_2 - B  x_1)^2 + \frac{ \lambda_1}{2} x_1^2 + \frac{ \lambda_2}{2} D_{\xi_2}^2 \,,\een

which after the dilation $\ds y_2=\frac{h\xi_2}{\sqrt{\lambda_2/2}}$ becomes
\ben  h^2D_{x_1}^2 + h^2D_{y_2}^2 + (\sqrt{\frac{ \lambda_2}{2}}y_2 - B  x_1)^2 + \frac{ \lambda_1}{2} x_1^2 \,.\een

A new change of coordinates leads to the sum of the two harmonic oscillators $h^2D_{z_j}^2+\mu_j z_j^2$, $j=1$, $2$, where
$\mu_1$, $\mu_2$ are the eigenvalues of the quadratic form
$(\sqrt{ \lambda_1/2}y_2 - B  x_1)^2 + \frac{ \lambda_1}{2} x_1^2$. These oscillators have discrete spectrum and their
lowest eigenvalues are $h\sqrt{\mu_j}$. A straightforward computation gives 
\ben (\sqrt{\mu_1}+\sqrt{\mu_2})^2=B^2+\lambda_{1,0}^2 \,. \een 
\qed \\
 
A result of  \cite{HeSjPise} allows then to estimate the first eigenvalue of a single well Schr\"odinger operator using the harmonic approximation. We also refer to \cite{cfks}, \S 11 for other results in this spirit.

  \begin{proposition} \label{harapp} 
Consider a vector field $\tilde A = (\tilde A_1, \tilde A_2) \in C^\infty (\R^2)$ and  a real nonnegative potential 
$ \tilde V \in C^{\infty} (\R^2)$  with a unique non degenerate minimum at a point $m \in \R^2$.  The smallest eigenvalue 
$\lambda_{ h,B } $ of the magnetic Schr\"odinger  operator
 \be \label{harmagapprox2}
	 P_{ h,\tilde A, \tilde V} = (hD_{x_1} - \tilde A_1(x))^2+ (hD_{x_2} - \tilde A_2(x))^2 + \tilde V(x)  
\ee 
 
 is simple and satisfies 
\ben
	 | \lambda_{ h,B }  -h \lambda_{ har, 1,B } |\leq Ch^{\nicefrac32} \,.
\een
Moreover, there exists $\epsilon_0>0$ such that $\sigma(P_{ h, \tilde A, \tilde V} ) \cap [0, h (\lambda_{ har, 1,B }+\epsilon_0)] = \{  \lambda_{ h,B }  \}$. 
\end{proposition}

 \begin{remark}
 In the case of a weak constant magnetic field $B=hB_0$, the harmonic approximation has no magnetic contribution and we have
 \ben
	 | \lambda_{ h,B }  -h \lambda_{ har, 1,0 } |\leq Ch^{\nicefrac32} \,.
\een
 \end{remark}

\subsection{Agmon distance}

Consider the Agmon metric $ V \, dx^2 $. For a piecewise $C^1$ curve $ \eta$, we can define its length $| \eta|$ in this metric,
 and for $x,y \in \R^2$ we define the Agmon distance $d_V(x,y)$ as the $\inf | \eta|$  over all  piecewise $C^1$ curves 
$ \eta$ joining $x$ to $y$. This distance may be degenerate in the sense that $d_V(x,y)=0$ for $x\neq y$, but it satisfies 
standard properties such as
\ben
	d_V(x,y) = d_V(y,x) \quad \text{ and } \quad d_V(x,z)\leq d_V(x,y)+ d_V(y,z) \,.
\een  

In the following, for $\varphi \in L^2(\R^2)$ and $y \in \R^2 $ we will use the notation 
\ben
	\varphi =   \mathcal{O}_\epsilon \left( e^{ \frac{-d_V (\cdot,y )(1-\epsilon) + \epsilon}h  } \right) \,,
\een	

which means that  for every $\epsilon>0$, there exists $h_\epsilon >0$ and $C_\epsilon$ such that 
\ben
	\left\| e^{ \frac{ d_V (\cdot,y )(1-\epsilon)  }h  }  \varphi (\cdot) \right\| _{L^2(\R^2)} \leq C_\epsilon e^{\frac \epsilon h}
\een

for   $h \in (0,h_\epsilon) $. Here $d_V (\cdot,m )$ is the Agmon distance to the point $m$. We refer to  \cite{DiSj}, \S 6
for details on Agmon estimates.

\subsection{Construction of a basis of the space attached to the  low lying spectrum of $P_{h,  A,V}$}

We now explain Carlsson's construction of an orthonormal basis of the spectral space attached to the low lying spectrum of $P_{h,  A,V}$ and prove Theorem \ref{theorem1}. The approach of Carlsson is quite general (no assumption of periodicity is needed) but he does not consider  the case with magnetic field. Nevertheless, the theory is simpler in the periodic case and it was shown in Section 9 of \cite{HeSjHarper1}  how to generalize this result with the help of \cite{HeSjPise}.\\
 
We follow the method of  ``filling the wells'' to obtain a basis of the spectral space attached to the low lying spectrum of $P_{h,  A,V}$. In our setting (see (\ref{minimaV})), the wells correspond to the points of the kagome lattice. The method consists then in associating with each $m \in \Gamma$, the Schr\"odinger operator $P_m$ given in (\ref{defPm}), which is obtained by filling all the other wells. Then, we get the desired basis considering the space spanned by the ground states of the $P_m$. Moreover, this basis respects the action of the magnetic operators, which lead to properties (\ref{Wtrans}) and (\ref{Wrot}). \\

 \textbf{Proof of Theorem \ref{theorem1}.} (\textit{Step 1})  Consider the operators $P_m$ defined in (\ref{defPm}). We have seen in Proposition \ref{propGG} that for all $m,n \in \Gamma$  there is $g  \in \mathcal{G}$ such that $g (m)=n$. Considering the associated $T_g$ defined in (\ref{worm}), all the operators $P_m$ are unitary equivalent.   \\

A result of Persson (\cite{persson}) gives that $\sigma(P_m )$ is discrete in the interval $[0,b]$ where $b$ is defined in
(\ref{definitionb}). Each operator $P_m$ is a Schr\"odinger operator with electric potential $V+V_m$. Using Hypothesis
\ref{hypothesisPotV},  $V+V_m$ has a unique non degenerate minimum, so Proposition \ref{harapp} applies to $P_m$.
The first eigenvalue $\lambda_{h,B}$ of  $P_m$ is simple and there exists $\epsilon_0>0$ such that  $\sigma(P_m) \cap I(h)  =  \{ \lambda_{h,B} \}$, where $I(h) =    [0, h (\lambda_{har,1,B(m) }+\epsilon_0)]$.\\

(\textit{Step 2})  Consider  $m_1 = m_{(0,0),1} $ and let $\varphi_1$ be an eigenfunction of $P_{m_1}$ with eigenvalue
 $\lambda(h)$ such that 
\be \label{choicevp1}
	\|  \varphi_1\|_{L^2(\R^2)} =1 \quad \text{ and } \quad  \varphi_1(m_1 ) \text{  is real.} 
\ee

For $\l = 1,3,5 $ we define
\be
	\varphi_\l  = F^{1-\l } \varphi_1 \,,
\ee

and  for every $  m_{\a,\l} \in    \Gamma$ we define an eigenfunction of $P_{  m_{\a,\l}} $ with eigenvalue $\lambda(h)$, by
\be \label{defvarphi}
	\varphi_{\tilde m_{\a,\l}}  = e^{-i\frac{\gamma}{2} \a \wedge \tilde \nu_\ell} \, T^\alpha \varphi_\l \,.
\ee

Defining $r_{ \tilde m} = (P-\lambda(h))\varphi_{ \tilde m}$, we have the Agmon estimates
\be \label{expdecay1}
	\varphi_{\tilde m } , r_{\tilde m} , \nabla_A\varphi_{\tilde m }  ,  \nabla_A  r_{\tilde m  }  = \mathcal{O}_\epsilon \left( e^{ \frac{-d_V (x,m )(1-\epsilon) + \epsilon}h  }  \right) \,,
 \ee 

where  $\nabla_A   =  (h D_{x_j}-A_j(x))_{j=1,2 }$.  Moreover,  
\be \label{expdecay2}
		\text{supp} \, r_{ \tilde m } \subset \bigcup_{n \in \Gamma \setminus\{m\} } B(n, \delta) \,. 
\ee 

 We also observe   by the harmonic approximation that
\be \label{hundemi}
	\varphi_1(\tilde m_1 ) =h^{-\nicefrac12} c_{har}  + \mathcal{O}(1) \,.
\ee

We now give the action of the magnetic operators over the eigenfunctions $\varphi_{ \tilde m}$.  The proof of the following Proposition is given at the end of this section.

\begin{proposition} \label{action2} For every $h>0$ there exist $c \in \{-1, 1\}$  such that for all $  m \in   \Gamma$ and
 $\b \in \Z^2$ we have
\ba
	T^\beta \varphi_{ \tilde m}  &=& e^{i\frac{\gamma}{2} \beta \wedge \tilde m} \, \varphi_{\tilde m+\b}  \,, \label{action2a} \\
	F \varphi_{ \tilde m }   &=&   c \, \varphi_{\kappa^{-1}(\tilde m) }  \label{action2b}   \,.
\ea
	  
\end{proposition}


(\textit{Step 3}) We may now state Carlsson's result. Let $\Sigma$ be  the spectral space associated with  $I(h)$ and $\Pi$ the
orthogonal projection over $\Sigma$. We define the projections 
\be \label{defvpro}
	v_{\tilde m}  =  \Pi \, \varphi_{\tilde m}   \,, \qquad  {  m}  \in   \Gamma \,.
\ee

 By estimates (\ref{expdecay1}) and (\ref{expdecay2}), for every $\epsilon>0 $ we can choose $\delta>0$ in the definition of
$P_m$ in (\ref{defPm}) such that 
\ban
 	| \left<v_{\tilde m} ,v_{ \tilde n} \right>| &\leq&  \exp \left(-\frac{(1-\epsilon) d_V(m,n)}h \right) \,,\\
| \left<v_{\tilde m} ,v_{ \tilde m} \right>-1| &\leq&  \exp \left(-\frac{(2S-\epsilon)}h \right)
 \ean
 
  for $h \in (0,h(\epsilon))$.  We denote  $D$ the matrix given by 
$D_{\tilde m,\tilde n}= \left<v_{\tilde n},v_{\tilde m} \right>$ and we define the functions 
\be \label{defesum} 
	e_{\tilde m} =  \sum_{ \tilde m \in \tilde \Gamma}  v_{\tilde n}(D^{-1/2})_{\tilde m,\tilde n} \,.
\ee 

The functions $e_{\tilde m} $ form an orthonormal basis of $\Sigma$.\\
Let $W_\gamma$ be the matrix of $P_{h,A,V}|\Sigma$ in this basis and put
 \ben
 	  (\tilde W)_{\gamma,\tilde m, \tilde n}=\left<r_{\tilde n},\varphi_{\tilde m}\right> \,.
\een

After Carlsson's theorem, for every $\epsilon>0 $ we can choose $\delta$ in the definition of $\chi$ in Step 1 such that
for $h \in (0,h(\epsilon))$
\be
\left|\left(W_\gamma\right)_{\tilde m,\tilde n} - \tilde W_{\gamma,\tilde m,\tilde n}\right| 
\leq  \exp \left(-\frac{(1-\epsilon)\, d^{(2)}_V( m, n)}h\right)
 \ee
 where
\ben d^{(2)}_V(n,m)=\min\{d_V(n,p)+d_V(p,m);\, p\in\Gamma,\,p\not=n,\,p\not=m\}.  \een
  
This proves (\ref{estcoeffintmatrix}) and (\ref{coefdiag}) in Theorem 1.3. \\
  
Moreover, the following proposition (which proof is given at the end of this section), proves that the orthonormalization process preserves the action of the magnetic operators. 
  
  \begin{proposition} \label{action3}
For every $  m  \in   \Gamma$ and $\beta \in \Z^2$ we have
\ba
	T^\b e_{ \tilde m}  &=&e^{i\frac{\gamma}{2} \b \wedge \tilde m} \, e_{\tilde m+\b} \label{action3a} \,,\\
	F e_{ \tilde m}  &=&   c \, e_{k^{-1}(\tilde m)} \label{action3b}   \,,
\ea
where $c  \in \{-1, 1\}$ is  defined in Proposition \ref{action2}.	
	  
\end{proposition}
  
Finally, properties (\ref{Wtrans}) and  (\ref{Wrot}) in Theorem \ref{theorem1} follow from Lemma \ref{dosh}, together with (\ref{action3a}) and (\ref{action3b}).  \\ \qed \\

 We end this section with the proofs of Propositions \ref{action2} and Proposition \ref{action3}.\\

\textbf{Proof of Proposition \ref{action2}.}  Let $  m=   m_{\a,\l}$. For the first relation, using (\ref{Talphaplusbeta}) we have
\ban
	T^\beta \varphi_{ \tilde m} &=& e^{-i\frac{\gamma}{2} \a \wedge \tilde \nu_\ell}  \, T^\beta T^\alpha \varphi_\l  \\
					  &=&  e^{-i\frac{\gamma}{2} \left(  \a \wedge \tilde \nu_\ell-\beta \wedge \alpha \right)}  \, T^{\alpha+\beta} \varphi_\l   \\
					  &=&  e^{i\frac{\gamma}{2}     \beta \wedge \tilde m} \, e^{-i\frac{\gamma}{2}    ( \a + \b)  \wedge \tilde \nu_\ell } \, T^{\alpha+\beta} \varphi_\l    \\
					  &=&   e^{i\frac{\gamma}{2} \beta \wedge \tilde  m} \,  \varphi_{\tilde m+\b} \,.
\ean

After (\ref{deftildenu}) and (\ref{relationsskappa}) we have 
\ban
	\kappa^{-1}(m_{\a,\l}) &=& \kappa^{-1}(\a) + \frac12 \kappa^{\l-2}(1,0)\\
				 &=& \kappa^{-1}(\a + \kappa^{\l-1}(1,0)) + \frac12 \kappa^{\l+1}( 1,0) \\
				 &=& \tilde m_{ \kappa^{-1}(\a + \kappa^{\l-1}(1,0)) ,\l+2} 
\ean

and
 \ben
	\kappa^{-1}(\a + \kappa^{\l-1}(1,0))   \wedge  \tilde \nu_{\l+2} =  - \kappa^{-1}(\a )   \wedge \frac12 \kappa^{\l-2}(1,0)  \,, 
\een

so we have to prove  that
\be \label{toproveFvarphi}
	F \varphi_{\tilde m} = c  \, e^{ i\frac{\gamma}{2}   \kappa^{-1}(\a )   \wedge \frac12 \kappa^{\l-2}(1,0)     } T^{\kappa^{-1}(\a + \kappa^{\l-1}(1,0)) } \,  \varphi_{\l+2}
\ee

for some $c \in \{ -1,1\}$. \\

Using (\ref{rtaualpha}) we have $t_1 r^{3}  V_{m_1} = V_{m_1}$ so $T_1F^3$ commutes with the multiplication by $V_{m_1}$.
 Lemma \ref{dosh} yields then that $T_1F^3$ commutes with $P_{m_1}$. Hence, since $\lambda(h)$ is a simple eigenvalue, there is
 a complex number $c$, such that  $|c|=1$ and
\be \label{pan}
	T_1F^3 	\varphi_1 = c \, \varphi_1 \,.
\ee

Moreover,  (\ref{relaF6}) and (\ref{relaFT}) yield $(T_1F^3)^2 = Id_{L^2(\R^2)}$, so $c^2=1$. Using (\ref{TalphaF}) and 
(\ref{pan}) we have
\ben
	F\varphi_1  = c  \, F^4  T^{-1}_1\varphi_1 =  c  \,  T^{\kappa^{-4}(-1,0)} F^{4}  \varphi_1 = c  \,  T^{\kappa^{-1}( 1,0)} F^{-2}  \varphi_1 \,.
\een

Considering (\ref{deftildenu}), (\ref{relationsskappa}), (\ref{Talphaplusbeta}) and (\ref{TalphaF}), the previous equality gives
\ban
	 F \varphi_{ \tilde m}   &=&  e^{-i\frac{\gamma}{2} \a \wedge \tilde \nu_\ell} \,  T^{\kappa^{-1}(\a)} F^{1-\l} F\varphi_1 \\
				&=& c \,   e^{-i\frac{\gamma}{2} \a \wedge \tilde \nu_\ell} \,  T^{\kappa^{-1}(\a)} F^{1-\l} T^{\kappa^{-1}( 1,0)} F^{-2}  \varphi_1 \\
				&=& c \,   e^{-i\frac{\gamma}{2} \a \wedge \tilde \nu_\ell} \,   T^{\kappa^{-1}(\a)}  T^{\kappa^{\l-2}(1,0)}_1  F^{1-(\l+2)}  \varphi_1 \\
				&=& c  \, e^{-i\frac{\gamma}{2}  \kappa^{-1}(\a)   \wedge  \left( \frac12 \kappa^{\l-2}(1,0) - \kappa^{\l-2}(1,0) \right)}   \, T^{\kappa^{-1}(\a + \kappa^\l(1,-1)) } \,  \varphi_{\l+2}  \,,
\ean

which yields (\ref{toproveFvarphi}) and ends the proof. \\ \qed\\

\textbf{Proof of Proposition \ref{action3}.} After Lemma \ref{dosh}, $T^\beta $ commutes with $P_{h,A  ,V}$, so also with $\Pi$ 
using the functional calculus of $P_{h,A  ,V}$. Then, using (\ref{action2a}), (\ref{defvpro}) and (\ref{defesum}) we get   
 \begin{eqnarray}
			T^\beta e_{ \tilde m}  &=&  \sum_{\tilde n}  T^\b v_{ \tilde n}\, (D^{-1/2})_{\tilde m, \tilde n} \nonumber \\
						        &=&  \sum_{\tilde n} e^{i\frac{\gamma}{2} \b \wedge \tilde n}\,  v_{\tilde n + \b } (D^{-1/2})_{\tilde m, \tilde n}   \label{eq:endsum} \\
						        &=& e^{i\frac{\gamma}{2}  \b \wedge \tilde m}  \,\sum_{\tilde n} v_{ \tilde n+\b}  \, e^{-i\frac{\gamma}{2}  \b \wedge \tilde m} \, (D^{-1/2})_{\tilde m, \tilde n} \, e^{i\frac{\gamma}{2} \b \wedge \tilde n} \nonumber \,. 
\end{eqnarray}

Since $T^\b$ is unitary, (\ref{action2a}) yields
 \begin{equation*}
\hat{D}_{\tilde m, \tilde n} :=\left< v_{ \tilde n+\b}, v_{ \tilde m+\b} \right>  =  e^{-i\frac{\gamma}{2}  \b \wedge \tilde m} \, D_{\tilde m  ,\tilde n}  e^{i\frac{\gamma}{2} \b \wedge \tilde n} \,. 
 \end{equation*}

Considering the diagonal matrix $A_{\tilde  m, \tilde m}  =   e^{-i\frac{\gamma}{2}  \b \wedge \tilde m} $, we note that
 \begin{equation*}
 	\left( (ADA^{-1})^{-1/2} \right)_{\tilde m,\tilde n} = 
 e^{-i\frac{\gamma}{2}  \b \wedge \tilde m} \, (D^{-1/2})_{\tilde m,\tilde n} \,  e^{i\frac{\gamma}{2} \b \wedge \tilde n} \,, 
 \end{equation*} 

so (\ref{eq:endsum}) becomes
\begin{eqnarray*}
			T^\b e_{ \tilde m}  &=& e^{i\frac{\gamma}{2}  \b \wedge \tilde m}  \,  \sum_{\tilde n}  v_{ \tilde n+\b} (\hat{D}^{-1/2})_{\tilde m,\tilde n} \,.\\ 
\end{eqnarray*}

We get (\ref{action3a}) noting that the sum in the right hand side of the previous equality is the $\tilde m+\b$ vector in the orthonormalization of $\{ v_{ \tilde n}\} $. \\

Similarly, using (\ref{action2b}) we find  
\be
	F e_{ \tilde m}  = c \, \sum_{\tilde n}  v_{\kappa^{-1}(\tilde  n)}  (D^{-1/2})_{\tilde m,\tilde n} \,,
\ee

and since the magnetic rotation is a unitary operator, we have $D=\tilde{D}$  where 
\be
	\tilde{D}_{\tilde m,\tilde n} =\left< v_{\kappa^{-1}(\tilde n)}, v_{\kappa^{-1}(\tilde m)} \right>\,.
\ee

Hence,
 \begin{eqnarray*}
	F e_{ \tilde m}  &=& c\, \sum_{\tilde n}  v_{\kappa^{-1}(\tilde n)} (\tilde{D}^{-1/2})_{\tilde m,\tilde n} \,.  
\end{eqnarray*}

We get (\ref{action3b}) reasoning as before.   \\ \qed

\section{The reduced models} \label{matrix}  

\subsection  {Introduction}  

In this section we obtain and study the reduced models associated with the low lying spectrum of $P_{h,A,V}$. First, under
Hypothesis \ref{hypothesisNN} and in the case of a weak and constant magnetic field, we estimate the coefficients of
$W_\gamma$ corresponding to the nearest neighbours for the Agmon distance. Then, by only keeping these terms, we construct
the operators $Q_{\gamma,\omega}$ and $P^{kay}_{\gamma,\omega}$ and prove Theorems \ref{studyhatW} and \ref{theopsudodiff}.
We then look at the case of rational values of the renormalized flux and reduce the  operator $W_\gamma$ to a family of
hermitian matrices. We end this article by proving the symmetries of the spectrum and the existence of eigenvalues and
flat bands.  

\subsection{The nearest neighbors and the tunneling effect}

As in \cite{HeSjHarper1} and \cite{kerdelhue1}, we implement here the results of \cite{HeSjPise} about the tunneling effect to estimate the coefficients of $W_\gamma$ corresponding to the nearest neighbours for the Agmon distance.\\

For any $\alpha \in \Z^2$ we want to identify in $W_\gamma$  the main terms corresponding to the interactions between the 
nearest wells for the Agmon distance to the triple $ \{m_{\a,1},m_{\a ,3},m_{\a ,5}\}$. After Hypothesis \ref{hypothesisNN} A, 
the nearest neighbours for the Agmon distance of a point $m_{\a ,j} \in \Gamma$ are (see Figure \ref{labeling}): 
\be 
	m_{( \a+ 2\tilde \nu_j),j-2}  \,,  m_{( \a- 2 \tilde \nu_{j-2}),j-2} \,, m_{( \a+ 2\tilde \nu_j),j+2}  \text{ and }  m_{( \a- 2\tilde \nu_{j+2}),j+2}   \,.  \quad \label{listeNN} 
\ee

\textbf{Proof of Theorems \ref{studyhatW} and \ref{theopsudodiff}}.  \textit{(Step 1)}  First, we notice that the term 
 $e^{-\frac{ i\gamma}{2}\alpha\wedge\tilde \nu_\ell}$ from the definition of the eigenfunctions in (\ref{defvarphi}) give 
the nice relations of Proposition \ref{action3}, but leads to the matrix $W_\gamma$ which does not satisfy the hypotheses 
of Theorem \ref{theopseudomatrix}.  To solve this, we introduce the new basis $\{f_{\tilde m} \}_{m \in \Gamma}$  
\be \label{defpsi}
	 f_{\tilde m_{\alpha,j}} = e^{i\frac{\gamma}{2}\alpha\wedge\tilde \nu_j} e_{\tilde m_{\alpha,j}} \,.
\ee
such that
\be \label{psitrans}
 T^\beta f_{\tilde m_{\alpha,j}}=e^{i\frac{\gamma}{2}\beta\wedge\alpha} f_{\tilde m_{\alpha+\beta,j}}\,.
\ee

Let $M$ be the matrix of $P_{h,A,V}-\lambda(h)\,I$ in this new basis. We denote by $M_{\alpha,\beta}^{j,k}$ the coefficients
 of $M$ in $\mathbb C$ and by $M_{\alpha,\beta}$ the blocks in $M_3(\mathbb C)$ 
(i.e.\ $\left(M_{\alpha,\beta}\right)^{j,k}=M_{\alpha,\beta}^{j,k}$). The matrix $M$ is obtained by conjugation by the diagonal
 matrix $A$ acting on  $\ell^2(\Z^2;\C^3) $ by
\ben
	A_{\alpha,\alpha}^{jj} = e^{i \frac \gamma2 \alpha \wedge \tilde \nu_j } \,, 
\een

so we obtain
\be \label{MtoW}\
	M_{\alpha,\beta}^{j,k}   = e^{-i\frac{\gamma}{2}(\alpha\wedge \tilde \nu_j-\beta\wedge \tilde \nu_k)} W_{\tilde m_{\alpha,j},\tilde m_{\beta,k}} \,.
\ee

The matrix $M$ thus inherits the properties of $W_\gamma$. Indeed, (\ref{Wtrans}) and (\ref{Wrot}) yield  
\ba  
	M_{\alpha,\beta} & = & e^{i\frac{\gamma}{2}(\alpha-\beta)\wedge\delta}M_{\alpha+\delta,\beta+\delta} \label{transM}  \,, \\
	M_{\alpha,\beta}^{j,k} & = &  e^{ i\frac{\gamma}{2}(  \kappa^{-1}(\alpha) \wedge \tilde \nu_{j+2} - \kappa^{-1}(\beta) \wedge \tilde \nu_{k+2})}    M^{j+2,k+2}_{\kappa^{-1}(\alpha)+\kappa^{j-2}(1,0),\beta+\kappa^{k-2}(1,0) } \label{rotM}  
\ea

for $\alpha$, $\beta$, $\delta\in\mathbb Z^2$ and $j$,$k\in\{1,3,5\}$, where $\kappa$ is defined in (\ref{defkappa}). \\

Relation (\ref{transM}) implies 
\ben 
	M_{\alpha,\beta}   =   e^{-i\frac\gamma2 \alpha\wedge\beta} M_{\alpha-\beta,0}  \,,
\een

which allows us to apply Theorem \ref{theopseudomatrix}. Hence, we define the operator 
\be \label{defkagrescomplet} 
	Q_\gamma =\sum_{\beta\in\mathbb Z^2} e^{i\frac{\gamma}{2}\beta_1\beta_2}\, M_{\beta,0}\,\tau_1^{\beta_1}\tau_2^{\beta_2}
 \ee

on $\ell^2(\Z^2;\C^3)$ and the symbol 
\be \label{defsymbkagomplet} 
	p(x,\xi,\gamma)=\sum_{\beta\in\mathbb Z^2} M_{\beta,0}\, e^{i(\beta_1x+\beta_2\xi)} \,,
\ee

and obtain that $W_\gamma$ and $Q_\gamma$ are unitary equivalent, and that the Weyl quantization $P_\gamma$ of $p(x,\xi,\gamma)$ and $W_\gamma$
have the same spectrum.\\

\textit{(Step 2)} We now compute the relations between the terms of $M_{\beta,0}$ corresponding to the interactions between
 the nearest wells of $V$. We have to consider 12 terms,  which correspond to the neighbours given in (\ref{listeNN}) for 
$\alpha=(0,0)$ and $j=1,3,5$. Relations (\ref{MtoW}) and (\ref{transM}) yield
\be \label{relMW} 
	M_{(-1,0),(0,0)}^{1,3} =   e^{i \frac \gamma4}  \,W_{\tilde m_{(0,0),1},\tilde m_{(1,0),3}} \,.
\ee 

Combining (\ref{rotM}) and (\ref{transM}) we find $M_{(-1,0),(0,0)}^{1,3} =  M_{(0,-1),(0,0)}^{3,5} $. Hence, using 
(\ref{rotM}) twice,   
\ben 
	M_{(-1,1),(0,0)}^{1,3} =  M_{(1,0),(0,0)}^{5,1}  =  M_{(0,-1),(0,0)}^{3,5}  = M_{(1,-1),(0,0)}^{3,5}  =  M_{(0,1),(0,0)}^{5,1}  =M_{(-1,0),(0,0)}^{1,3}   \,.
\een

By the self-adjointness of $P_{h,A,V}$, the other six terms equal the complex conjugate of $M_{(-1,0),(0,0)}^{1,3}$.  \\

\textit{(Step 3)} We now estimate $M_{(-1,0),(0,0)}^{1,3} $ when $A$ is given by (\ref{weakA}) and $V$ satisfies
Hypothesis~\ref{hypothesisNN}. First,
\be \label{fluxmagcte}
	\gamma     =  B_0\,  (2\nu_1) \wedge (2\nu_2)  = 2 \sqrt 3 B_0 \,. 
 \ee

We compute explicitly the value of the phases of the magnetic translations and rotation. Remarks \ref{remark1} and 
\ref{remark2}, together with (\ref{fluxmagcte}), give that for any  $\alpha \in \Z^2$ and $\varphi \in L^2(\R^2)$ 
 \be \label{TFjaar}	
 	(T^\alpha \varphi)(x) =  e^{ -i \frac {B_0}2 \, x \wedge \left( 2\alpha_1\nu_1 +2\alpha_2\nu_2 \right)   } \, \varphi\left( \tau^{-\alpha}(x) \right) \quad \text{ and } \quad (F\varphi)(x) = \varphi(r(x)) \,.
 \ee

 The results in Section 3 of \cite{HeSjPise} give an asymptotic estimate for 
$(W_\gamma)_{\tilde m_{(0,0),1} , \tilde m_{(1,0),3}} $. In order to apply the results there,  we first need to verify that 
the values of the functions $\varphi_{ \tilde m_{(0,0),1}}$ and $\varphi_{\tilde m_{(1,0),3}}$ at the bottom of their 
respective wells are real. The value of $\varphi_{\tilde m_{(0,0), 1}} (m_{(0,0),1})= \varphi_1 ( \nu_1) $ has been chosen 
real in (\ref{choicevp1}). The definition  of  $\varphi_{ \tilde m}$ in (\ref{defvarphi}) and $m_{(1,0),3} =  \nu_1+\nu_2 $, 
together with (\ref{fluxmagcte}) and (\ref{TFjaar}), give   
\ban
	\varphi_{\tilde m_{(1,0),3}}  (m_{(1,0),3}) &=& e^{-i \frac{\gamma}{2 }\,  2 \tilde \nu_1 \wedge \tilde \nu_3 } \,  \left( T_1F^{-2} 	\varphi_1 \right) (\nu_2) \\
						&=& e^{ - i \frac  {B_0}4  \,  (2\nu_1) \wedge (2\nu_2) } \, e^{-i  \frac{B_0}2 \, \nu_2 \wedge 2\nu_1  } \,     \varphi_1    ( r^{-2}(\nu_3) ) \\
						&=&     \varphi_1 ( \nu_1)    \,,
\ean

so $\varphi_{\tilde m_{(1,0),3}} (m_{(1,0),3}) $ is also real. \\

The results  in \cite{HeSjPise} are given for a magnetic potential $A_t = t A$, with the condition (see (2.40) therein) 
$ |t| = \mathcal{O}(h^{\nicefrac12} (-\ln h)^{\nicefrac12}  ) $. Our setting satisfies this requirement with $t = h$.
 By Proposition 3.12, Remark 3.17 and Lemma 3.15 in \cite{HeSjPise} \footnote{Formula (3.26) in \cite{HeSjPise} has 
unfortunately disappeared in the printing and reads: $W_{jk}^t = h^{\nicefrac12 } b_{jk}^t (h) e^{-S_{jk}/h} $ (\cite{Hpc}).} 
 and  assuming Hypothesis~\ref{hypothesisNN} we get 
\ben
	 (W_\gamma)_{\tilde m_{(0,0),1}, \tilde m_{(1,0),3}} = h^{\nicefrac 12} \, b(h) \, e^{- \frac {S(h)}h} \,,
\een 

where  
 \ben  
 \begin{array}{rl} 
	|b(h)| =  b_0 + \mathcal{O}(h) \,, &\text{Re}(S(h) )=  d_V(m_{(0,0),1},m_{(1,0),3}) + \mathcal{O}(h^2)  \,,\\
	\text{Arg}(b(h) ) = \, \, \pi+ \mathcal{O}(h)  \,, & \text{Im}(S(h) ) =  \text{Circ}(A,\zeta)+ \mathcal{O}(h^3)  \,.
 \end{array} 
 \een

Here before $\zeta : [0,1] \to [m_{(0,0),1}, m_{(1,0),3} ]$ is the unique minimal geodesic  between $m_{(0,0),1}$ and
 $ m_{(1,0),3}$ and 
\ben
	\text{Circ}(A,\zeta) = \int_\zeta  \omega_A \,.
\een

Considering (\ref{weakA}) and (\ref{fluxmagcte}) we obtain $\text{Circ}(A,\zeta) =  h \sqrt 3   B_0/4 = h \gamma/8$, so  
\ben
	 (W_\gamma)_{\tilde m_{(0,0),1}, \tilde m_{(1,0),3}} =   h^{\nicefrac 12} \, b_0 \, e^{- \frac Sh} \, (1+\mathcal{O}(h))  \, e^{i\left(-\frac \gamma 8+ \pi  + \mathcal{O}(h)\right)} \,.
\een 
 
After (\ref{relMW})  we find 
\be \label{valueM13} 
	M_{(-1,0),(0,0)}^{1,3} =  - \rho e^{i \omega}   e^{i\frac \gamma 8}   
\ee 

with $\rho $ and $\omega$ satisfying (\ref{estrhoh}) and (\ref{estbetah}).\\

\textit{(Step 4)} Equality (\ref{valueM13}) gives that the operator obtained when only considering the 12 terms from Step 2 in
the sum (\ref{defkagrescomplet}) equals $- \rho\, Q_{\gamma,\omega}$. Similarly, the symbol obtained when only considering in
  (\ref{defsymbkagomplet}) these terms equals $p^{kag}(x,\xi,\gamma,\omega)$. Finally, defining 
\ben
	-\rho\, R_\gamma = Q_\gamma -  (\lambda(h) I - \rho Q_{\gamma,\omega} ) 
\een	

the estimations (\ref{estcoeffintmatrix}) and (\ref{coefdiag}) in Theorem \ref{theorem1} achieve the proof.

\qed\\


{\bf Proof of Proposition \ref{gammadependence}} We first observe that the magnetic translations $\tau_1$ and $\tau_2$ do not change if we add to $\gamma$ a multiple of $2\pi$.
Hence, the spectrum of $Q_{\gamma+6\pi,\omega+\frac{\pi}4}$ is that one of
\ben \left(\begin{array}{ccc}
       0 & -e^{i(\omega+\frac\gamma8)} \left(\tau_1^*-e^{-i\frac\gamma2}\tau_1^*\tau_2\right) & -e^{-i(\omega+\frac\gamma8)} \left(\tau_1^*+\tau_2^*\right)\\
-e^{-i(\omega+\frac\gamma8)} \left(\tau_1-e^{-i\frac\gamma2} \tau_1\tau_2^*\right) & 0 & -e^{i(\omega+\frac\gamma8)} \left(\tau_2^*-e^{-i\frac\gamma2}\tau_1\tau_2^*\right)\\
     -e^{i(\omega+\frac\gamma8)} \left(\tau_1+\tau_2\right) & -e^{-i(\omega+\frac\gamma8)} \left(-e^{-i\frac\gamma2}\tau_1^*\tau_2+\tau_2\right) & 0
\end{array}\right)
\een
acting on $\ell^2(\mathbb Z^2;\,\mathbb C^3)$.\\
Using Theorem \ref{theopseudomatrix}, it is also the spectrum of the $\gamma$-quantized of the symbol
 \ben 
\left(\begin{array}{ccc}
       0 & e^{i(\omega+\frac\gamma8)} \left(-e^{-ix}+e^{-i(x-\xi)}\right) & e^{-i(\omega+\frac\gamma8)} \left(-e^{-ix}-e^{-i\xi}\right)\\
e^{-i(\omega+i\frac\gamma8)} \left(-e^{ix}+e^{i(x-\xi)}\right) & 0 & e^{i(\omega+\frac\gamma8)} \left(e^{i(x-\xi)}-e^{-i\xi}\right)\\
     e^{i(\omega+\frac\gamma8)} \left(-e^{ix}-e^{i\xi}\right) & e^{-i(\omega+\frac\gamma8)} \left(e^{-i(x-\xi)}-e^{i\xi}\right) & 0
\end{array}\right) 
\een
We simply compose this symbol by the affine symplectic map $(x,\xi)\mapsto (x+\pi,\xi+\pi)$ to recover
$p^{kag}_{\gamma,\omega}$.
At the level of pseudodifferential operators, this composition is associated with the conjugation by
$\ds u(\cdot)\mapsto e^{\frac{i}{\gamma}\pi\cdot}u(\cdot-\pi)$. This yields (\ref{gammatrans}).\\
We recall the operator $\Gamma\,:\,u(\cdot)\mapsto\bar u(-\cdot)$ introduced in the proof of (\ref{symhex}).
If $a$ is a symbol,
\ben \Gamma\,({\rm Op}^W_\gamma a(x,\xi))\,\Gamma={\rm Op}^W_\gamma \bar a(-x,\xi)\,. \een
This, together with the observation
\ben P^{kag}_{-\gamma,-\omega}={\rm Op}^W_\gamma p^{kag}(x,-\xi,-\gamma,-\omega)
={\rm Op}^W_\gamma \bar p^{kag}(-x,\xi,\gamma,\omega)\een
yield (\ref{gammasym}).

\subsection{Study of the spectrum for rational values of the renormalized magnetic flux}

 We now prove Theorem \ref{theoremmatsym}, which will allow us to numerically compute the spectrum of $Q_{\gamma,\omega}$. \\
 
\textbf{Proof of Theorem \ref{theoremmatsym}.} Using the definitions of $\tau_1$ and $\tau_2$ in (\ref{def1tau12}), we explicitly
write $Q_{\gamma,\omega}$ as
\ban 
	(Q_{\gamma,\omega} v)_\a^1 &=&    e^{i(\omega+\frac\gamma8)}\left(   v_{\a_1+1,\a_2}^3  +  e^{ - i \frac{\gamma}2 } e^{  i  \gamma(  \a_1+1)  }  \, v_{\a_1+1 ,\a_2-1}^3 \right)\\
&&  +  e^{-i(\omega+\frac\gamma8)} \left(   v_{\a_1-1,\a_2}^5 +  e^{ - i  \gamma   \a_1   }   \,  v_{\a_1,\a_2+1}^5   \right)   \\
	(Q_{\gamma,\omega} v)_\a^3 &=&    e^{-i(\omega+\frac\gamma8)}\left(   v_{\a_1+1,\a_2}^1  + e^{ - i \frac{\gamma}2 } e^{  -i  \gamma(  \a_1-1)  }  \, v_{\a_1-1 ,\a_2+1}^1  \right)  \\
	&& +  \,  e^{i(\omega+\frac\gamma8)}\left(   e^{ - i \frac{\gamma}2 } e^{  -i  \gamma(  \a_1-1)  }  \, v_{\a_1-1 ,\a_2+1}^5 + ^{ - i  \gamma   \a_1  }    v_{\a_1 ,\a_2+1}^5     \right)  \\
	(Q_{\gamma,\omega} v)_\a^5 &=&  e^{i (\omega+\frac\gamma8)} \left(    v_{\a_1-1,\a_2 }^1 +    e^{   i  \gamma  \a_1  } \,   v_{\a_1,\a_2-1}^1  \right) \\
&& +  e^{-i(\omega+\frac\gamma8)} \left(   e^{ - i \frac{\gamma}2 } e^{  i  \gamma(  \a_1+1)  }  \, v_{\a_1+1 ,\a_2-1}^5 +  e^{    i  \gamma  \a_1  } \,   v_{\a_1,\a_2-1}^5    \right) \,.  
\ean

We first notice that $Q_{\gamma,\omega} $ commutes with the translation $v_\a \mapsto v_{\a_1,\a_2-1}$, so we may use a Floquet
theory to obtain
\ben
	 \sigma(Q_{\gamma,\omega}) =  \bigcup_{ \theta_2 \in [0,1[  }  \sigma\left(  Q_{ \gamma,\omega, \theta_2}  \right)   \,,
\een 

where $Q_{ \gamma,\omega, \theta_2}  v = Q_{\gamma,\omega}  v  $ and 
\ben
	\mathcal{D}(Q_{ \gamma,\omega, \theta_2} ) = \left\{ v  : \Z^2 \to  \C^3 \in   \ell^2(  \Z_{\a_1} ;  \C^3)  \,;\, v_{\a_1,\a_2-1}   = e^{i2\pi\theta_2} v_\a \right\} \,.
\een	

 Since any sequence in $\mathcal{D}(Q_{ \gamma,\omega, \theta_2} ) $ is only determined by the first coordinate $\alpha_1$,
the operator $Q_{ \gamma,\omega, \theta_2}$ has the same spectrum that the operator $\check Q_{\gamma,\omega,\theta_2} $
acting on $\ell^2(\Z;\C^3)$ by 
\ban 
	(\check Q_{\gamma,\omega} v)_\a^1 &=&    e^{i(\omega+\frac\gamma8)}\left(   v_{\a-1 }^3  +  e^{ - i \frac{\gamma}2 } e^{  i  \gamma(  \a+1)  }  \, e^{i2\pi\theta_2} v_{\a+1 }^3 \right) \\
&& +  e^{-i(\omega+\frac\gamma8)} \left(   v_{\a-1}^5 +  e^{ - i  \gamma   \a   }  e^{-i2\pi\theta_2} \,  v_{\a }^5   \right)   \\
	(\check Q_{\gamma,\omega} v)_\a^3 &=&    e^{-i(\omega+\frac\gamma8)}\left(   v_{\a+1 }^1  + e^{ - i \frac{\gamma}2 } e^{  -i  \gamma(  \a-1)  } e^{-i2\pi\theta_2} \, v_{\a-1 }^1  \right)  \\
	&& +  \,  e^{i(\omega+\frac\gamma8)}\left(   e^{ - i \frac{\gamma}2 } e^{  -i  \gamma(  \a-1)  } e^{-i2\pi\theta_2} \, v_{\a-1 }^5 + ^{ - i  \gamma   \a  }  e^{-i2\pi\theta_2} \,  v_{\a }^5     \right)  \\
	(\check Q_{\gamma,\omega} v)_\a^5 &=&  e^{i (\omega+\frac\gamma8)} \left(    v_{\a-1 } +    e^{   i  \gamma  \a  } e^{i2\pi\theta_2} \,   v_{\a }^1  \right) \\
&& +  e^{-i(\omega+\frac\gamma8)} \left(   e^{ - i \frac{\gamma}2 } e^{  i  \gamma(  \a+1)  }  e^{i2\pi\theta_2}\, v_{\a+1 }^5 +  e^{  - i  \gamma  \a  } \,  e^{ i2\pi\theta_2}  v_{\a }^5    \right) \,.  
\ean

Now, if $\gamma / (2\pi) =   p/q$ then $\check Q_{\gamma,\omega,\theta_2}$ commutes with $\tau_1^q$.  We may then use another 
Floquet theory to obtain 
\be \label{sigmagamma123}
	 \sigma(Q_{\gamma,\omega}) =  \bigcup_{(\theta_1,\theta_2) \in [0,1[ \times [0,1[ }  
\sigma\left(  \check Q_{p,q,\omega,\theta_1,\theta_2}  \right)   
\ee 
 
where $\check Q_{p,q,\omega,\theta_1,\theta_2}  v = \check Q_{\gamma,\omega,\theta_2} v$ and 
 \ben
 	\mathcal{D}( \hat W_{p,q,\omega,\theta_1,\theta_2} ) = \left\{ v  \in \l^\infty(\Z ; \C^3)    \,;\, v_{\a+q}   = e^{i   2\pi  \theta_1 q } v_\a  \right\}\,.
 \een

Since any sequence  in $ \mathcal{D}(\check Q_{p,q,\omega,\theta_1,\theta_2} )$ is only determined by its $q$ first terms,
the operator $\check Q_{p,q,\omega,\theta_1,\theta_2}$ has the same spectrum that its restriction to $\C^{3\times q}$. Taking
in account the condition $v_{\a+q}   = e^{i   2\pi  \theta_1 q } v_\a$ the operator $\check Q_{p,q,\omega,\theta_1,\theta_2}$ has
the same spectrum that the operator  $\hat M_{p,q,\omega,\theta_1,\theta_2}$ acting on    
\ben
	v=(v^1,v^3,v^5) = (v^1_0, \cdots, v^1_{q-1},v^3_0, \cdots, v^3_{q-1},v^5_0, \cdots, v^5_{q-1})
\een

by
\ben 
 \hat M_{p,q,\omega,\theta_1,\theta_2}  = \left(\begin{array}{c|c|c} 
	&   &   \\
	  0_q   &  \hat M_{p,q,\omega,\theta_1,\theta_2}^{13}   & \hat  M_{p,q,\omega,\theta_1,\theta_2}^{15}   \\
	  &   &   \\
	\hline 
	  &   &   \\
	 \left[  \hat M_{p,q,\omega,\theta_1,\theta_2}^{13} \right]^*   &   0_q   &\hat  M_{p,q,\omega,\theta_1,\theta_2}^{35}    \\
	  &   &   \\
	\hline 
	 &   &   \\
	\left[  \hat M_{p,q,\omega,\theta_1,\theta_2}^{15} \right]^* & \left[  \hat M_{p,q,\omega,\theta_1,\theta_2}^{35} \right]^*  &   0_q   \\
	  &   &   \\
	\end{array}\right) \,,
\een

where
\ba  
	\hat M_{p,q,\omega,\theta_1,\theta_2}^{13}   &=&  e^{i (\omega+\frac{\pi}{4}\frac pq)}  (  e^{i2\pi\theta_1} \tilde K_q + e^{ i\pi\frac pq} e^{i  2\pi (\theta_1+\theta_2)}    J_{p,q}  \tilde K_q) \nonumber\\
	\hat M_{p,q,\omega,\theta_1,\theta_2}^{15}   &=&  e^{-i(  \omega+\frac{\pi}{4}\frac pq)}  ( e^{i2\pi\theta_1} \tilde K_q  +e^{-i2\pi\theta_2} J^*_{p,q})  \qquad  \label{Mpqtheta12in} \\
	\hat M_{p,q,\omega,\theta_1,\theta_2}^{35}   &=&   e^{i( \omega+\frac{\pi}{4}\frac pq)}  ( e^{ i\pi\frac pq} \, e^{ -2i\pi (\theta_1+\theta_2)}  J^*_{p,q}  \tilde K_q^* +  e^{- i2\pi\theta_2} J^*_{p_q}    )   \nonumber
\ea 

 with  $J_{p,q} $ defined in (\ref{definitpnjandk}) and 
 \ben
	 (  \tilde K_q   )_{ij} = \left\{ \begin{array}{cl} 1 & \text{if }  j=i+1 \, (\text{mod } q) \text{ and }  i \neq q \\  
				e^{iq 2\pi \theta_1}	& \text{if }  j=i+1 \, (\text{mod } q) \text{ and }  i =q  \\
				0  & \text{if }  j\neq i+1 \, (\text{mod } q)   
					\end{array}\right. \,.
\een	  
 
Since $J_{p,q} $ and $\tilde K_q$ satisfy the commutation relation 
\ben
	J_{p,q} \tilde K_q =  e^{ - i 2\pi   \frac pq  }  \tilde K_q  J_{p,q}  \,, 
\een
 
 we may rewrite  
\ba  
	\hat M_{p,q,\omega,\theta_1,\theta_2}^{13}   &=&  e^{i (\omega+\frac{\pi}{4}\frac pq)}  (  e^{i2\pi\theta_1} \tilde K_q + e^{-i\pi\frac pq} e^{i  2\pi (\theta_1+\theta_2)}     \tilde K_q J_{p,q} ) \,, \nonumber\\
 	\hat M_{p,q,\omega,\theta_1,\theta_2}^{35}   &=&   e^{i( \omega+\frac{\pi}{4}\frac pq)}  ( e^{ -i\pi\frac pq} \, e^{ -2i\pi (\theta_1+\theta_2)}   \tilde K_q^* J^*_{p,q} +  e^{- i2\pi\theta_2} J^*_{p_q}    )  \,. \nonumber
\ea 
 
Finally, noting that $ e^{i  2\pi  \theta_2}   J_{p,q} $ and $\tilde K_q$ are respectively the conjugate of
$ e^{i  2\pi  \theta_2}   J_{p,q} $ and $ e^{i  2\pi  \theta_1}   K_q $ by  the unitary matrix
$ \text{diag}(\exp{(2i\pi \theta_1(j-1) )}) $,  we have that $\hat M_{p,q,\omega,\theta_1,\theta_2}$ is unitary equivalent 
with the matrix in  (\ref{Mpqthetaijin}), which yields the proof. \hfill \qed \\

 We end this article by proving the symmetries of $\sigma_{\gamma,\omega}$ and the existence of eigenvalues.\\ 
 
 \textbf{Proof of Proposition \ref{propsymmkag}.} 
$Q_{\g} $ is the sum of four unitary operators, so (\ref{kagrange}) holds. As in the proof of Proposition
(\ref{gammadependence}), we observe that the magnetic translations $\tau_1$ and $\tau_2$ do not change if we add  to $\gamma$ a multiple of $2\pi$, so $Q_{\gamma+8\pi,\omega}=-Q_{\gamma,\omega}$, which yields (\ref{kagreflex4}) and thus 
(\ref{kagtrans}) is a direct consequence of (\ref{kagreflex4}). Combining (\ref{kagtrans}) and (\ref{kagreflex4}) with
 Proposition \ref{gammadependence}, we easily obtain (\ref{kagreflex1}), (\ref{kagreflex5}), (\ref{kagreflex11}) 
and (\ref{kagreflex12}). \\ \qed \\

 \textbf{Proof of Proposition \ref{eigenvalue}.} Using Theorem \ref{theoremmatsym} the spectrum of $Q_{0,0}$ is the union over
$(\theta_1,\theta_2)$
running on $[0,1]\times[0,1]$ of the spectra of the matrices 
\ben M_{0,1,0,\theta_1,\theta_2}=\left(
\begin{array}{ccc}
 0& e^{i2\pi\theta_1}+e^{i2\pi(\theta_1+\theta_2)} & e^{i2\pi\theta_1}+e^{-i2\pi\theta_2} \\
e^{-i2\pi\theta_1}+e^{-i2\pi(\theta_1+\theta_2)} & 0 & e^{-i2\pi(\theta_1+\theta_2)}+e^{-i2\pi\theta_2} \\
e^{-i2\pi\theta_1}+e^{i2\pi\theta_2}  & e^{i2\pi(\theta_1+\theta_2)}+e^{i2\pi\theta_2} &0
\end{array}\right) \een
which characteristic polynomial are $(\lambda+2)((\lambda-1)^2-(3+2p^\triangle (2\pi\theta_1,-2\pi\theta_2)))$. 
Since the range of $p^\triangle$ is $[-3/2,3]$, the three eigenvalues of $M_{0,1,\theta_1,\theta_2}$ respectively run on
$\{-2\}$, $[-2,1]$ and $[1,4]$.\\

Similarly, the spectrum of $Q_{4\pi,0}$ is the union over
$(\theta_1,\theta_2)$
running on $[0,1]\times[0,1]$ of the spectra of the matrices 
$M_{2,1,0,\theta_1,\theta_2}$ given by
\ben \left(
\begin{array}{ccc}
 0& i(e^{i2\pi\theta_1}+e^{i2\pi(\theta_1+\theta_2)}) & -i(e^{i2\pi\theta_1}+e^{-i2\pi\theta_2}) \\
-i(e^{-i2\pi\theta_1}+e^{-i2\pi(\theta_1+\theta_2)}) & 0 & i(e^{-i2\pi(\theta_1+\theta_2)}+e^{-i2\pi\theta_2}) \\
i(e^{-i2\pi\theta_1}+e^{i2\pi\theta_2})  & -i(e^{i2\pi(\theta_1+\theta_2)}+e^{i2\pi\theta_2}) &0
\end{array}\right) \een
which characteristic polynomial are $\lambda(\lambda^2-(6+2p^\triangle (2\pi\theta_1,-2\pi\theta_2)))$. 
So the three eigenvalues of $M_{2,1,0,\theta_1,\theta_2}$ respectively run on
$[-2\sqrt{3},-\sqrt{3}]$, $\{0\}$ and $[\sqrt3,2\sqrt3]$.\\

We compute the other characteristic polynomial using the symbolic computation software Mathematica. We obtain 
\ben\det(\lambda\,Id-M_{2,3,0,\theta_1,\theta_2})=
 (\lambda+\sqrt{3})^3(\lambda^6-3\sqrt{3}\lambda^5+18\sqrt{3}\lambda^3-36\lambda^2+6
-2p^{\triangle}(6\pi\theta_1,-6\pi\theta_2)) \,,\een
\begin{multline*}\det(\lambda\,Id-M_{4,3,0,\theta_1,\theta_2})=\\
(\lambda+1)^3(\lambda^6-3\lambda^5-12\lambda^4+38\lambda^3+24\lambda^2-120\lambda +70
-2p^{\triangle}(6\pi\theta_1,-6\pi\theta_2)) \,,\end{multline*}
\begin{multline*}\det(\lambda\,Id-M_{1,2,\frac{\pi}{8},\theta_1,\theta_2})=\\
(\lambda+\sqrt{2})^2(\lambda^4-2\sqrt{2}\lambda^3-6\lambda^2+12\sqrt{2}\lambda-6
+2p^{\triangle}(4\pi\theta_1,-4\pi\theta_2)) \,, \end{multline*} 
\ban \det(\lambda\,Id-M_{3,2,\frac{\pi}{8},\theta_1,\theta_2})&=&
(\lambda+2)^2(\lambda^4-4\lambda^3+8\lambda-2
+2p^{\triangle}(4\pi\theta_1,-4\pi\theta_2)) \,, \\
\det(\lambda\,Id-M_{-1,6,\frac{\pi}{8},\theta_1,\theta_2})&=&
\left(\lambda+\frac{-\sqrt{2}+\sqrt{6}}{2}\right)^6 ( T(\lambda) + 2 p^\triangle(12\pi\theta_1,-12\pi\theta_2) )\,,\\
\det(\lambda\,Id-M_{7,6,\frac{\pi}{8},\theta_1,\theta_2})&=&
\left(\lambda+\frac{\sqrt{2}+\sqrt{6}}2\right)^6 (U(\lambda)+2p^\triangle (12\pi\theta_1,-12 \pi \theta_2))\,,\ean

with
\ban T(\lambda) = &-&9726-5616 \sqrt{3}+(9828 \sqrt{2}+5652 \sqrt{6}) \lambda+(3024+1836 \sqrt{3}) \lambda^2\\
&-&(8244 \sqrt{2}+4596 \sqrt{6}) \lambda^3
+(1584+720 \sqrt{3}) \lambda^4+(2970 \sqrt{2}+1350 \sqrt{6}) \lambda^5\\
&-&(828+540 \sqrt{3}) \lambda^6
-(612 \sqrt{2}+144 \sqrt{6}) \lambda^7
+(36+180 \sqrt{3}) \lambda^8\\
&+&(38 \sqrt{2}+18 \sqrt{6}) \lambda^9+(6-21 \sqrt{3}) \lambda^{10}
+(3 \sqrt{2}-3 \sqrt{6}) \lambda^{11}+\lambda^{12} \,, \ean
and
\ban 
U(\lambda)&=&  -9726 + 5616 \sqrt{3} +36 \sqrt{2} (-273+157\sqrt{3}) \, \lambda + 
 108 (28 - 17 \sqrt{3}) \, \lambda^2 \\
&&+ 12 \sqrt{2} (687 - 383 \sqrt{3}) \,  \lambda^3  +   144 (11 - 5 \sqrt{3}) \, \lambda^4 + 270 \sqrt{2} (-11 + 5 \sqrt{3}) \,
 \lambda^5 \\
&&+ 36 (-23 + 15 \sqrt{3}) \, \lambda^6  
+ 36 \sqrt{2} (17 - 4 \sqrt{3}) \, \lambda^7
  +  36 (1 - 5 \sqrt{3})\, \lambda^8  \\
&& +2 \sqrt{2} (-19 + 9 \sqrt{3}) \,\lambda^9 + 
 3 (2 + 7 \sqrt{3}) \,\lambda^{10}- 3 \sqrt{2}(1 + \sqrt{3}) \,  \lambda^{11} + \lambda^{12}   \,. 
\ean  \qed \\

\bibliographystyle{alpha}
\bibliography{bib}

\newcommand{\etalchar}[1]{$^{#1}$}
\begin{thebibliography}{DFE{\etalchar{+}}05}

\bibitem[AJ09]{AvJi}
A.~Avila and S.~Jitomirskaya.
\newblock {The Ten Martini Problem}.
\newblock {\em Ann. of Math.}, 170:303--342, 2009.

\bibitem[Bel87]{Be1}
J.~Bellissard.
\newblock C$^*$-algebras in solid state physics-2d electrons in a uniform
  magnetic field.
\newblock {\em Warwick conference on operator algebras}, 1987.

\bibitem[BKS91]{BellissardKreftSeiler}
J.~Bellissard, C.~Kreft, and R.~Seiler.
\newblock Analysis of the spectrum of a particle on a triangular lattice with
  two magnetic fluxes by algebraic and numerical methods.
\newblock {\em J.Phys. A}, (10):2329--2353, 1991.

\bibitem[BS82]{BeSi}
J.~Bellissard and B.~Simon.
\newblock {Cantor spectrum for the almost Mathieu equation}.
\newblock {\em J. Funct. Anal.}, 48(3):408--419, 1982.

\bibitem[Car65]{Cartier}
P.~Cartier.
\newblock Quantum mechanical commutation relations and $\theta$ functions.
\newblock {\em Proc. Symp. Pure Math, Boulder}, pages 183--186, 1965.

\bibitem[CFKS87]{cfks}
H.L. Cycon, R.~Froese, W.~Kirsch, and B.~Simon.
\newblock {\em Schr{\"o}dinger operators}.
\newblock Springer, 1987.

\bibitem[CW79]{ClaroWannier}
F.~H. Claro and G.~H. Wannier.
\newblock Magnetic subband structure of electrons in hexagonal lattices.
\newblock {\em Phys. Rev. B}, 19(12):6068--6074, 1979.

\bibitem[DFE{\etalchar{+}}05]{damski}
B.~Damski, H.~Fehrmann, H.-U. Everts, M.~Baranov, L.~Santos, and M.~Lewenstein.
\newblock Quantum gases in trimerized kagome lattices.
\newblock {\em Phys. Rev. A}, 72, 2005.

\bibitem[DGJO11]{DaGerJuOh}
J.~Dalibard, F.~Gerbier, G.~Juzeli\={u}nas, and P.~\"{O}hberg.
\newblock Colloquium: Artificial gauge potentials for neutral atoms.
\newblock {\em Rev. Mod. Phys.}, 83(4):1523--1543, 2011.

\bibitem[DS99]{DiSj}
M.~Dimassi and J.~Sj\"ostrand.
\newblock {\em Spectral Asymptotics in the Semi-Classical Limit}.
\newblock Cambridge University Press, 1999.

\bibitem[Gei11]{Geim}
A.~K. Geim.
\newblock {Nobel Lecture: Random walk to graphene}.
\newblock {\em Rev. Mod. Phys.}, 83(3):851--862, 2011.

\bibitem[HA10]{HuAlt}
S.~D. Huber and E.~Altman.
\newblock {B}ose condensation in flat bands.
\newblock {\em Phys. Rev. B}, 82(184502), 2010.

\bibitem[Hel]{Hpc}
B.~Helffer.
\newblock Personnal communication.

\bibitem[Hel09]{HelfferSyrie}
B.~Helffer.
\newblock {\em Introduction to semi-classical methods for the {S}chr\"odinger
  operator with magnetic fields. In Aspects th\'eoriques et appliqu\'es de
  quelques EDP issues de la g\'eom\'etrie ou de la physique}.
\newblock Proceedings of the CIMPA School held in Damas (Syrie). S\'eminaires
  et Congr\`es. SMF, 2009.

\bibitem[Hel13]{HeBook}
B.~Helffer.
\newblock {\em Spectral Theory and its Applications}, volume 139 of {\em
  Cambridge Studies in Advanced Mathematics}.
\newblock Cambridge University Press, 2013.

\bibitem[Hof76]{hofstadter}
D.~R. Hofstadter.
\newblock Energy levels and wave functions of {B}loch electrons in rational and
  irrational magnetic fields.
\newblock {\em Phys. Rev. B}, 14(6):2239--2249, 1976.

\bibitem[Hou09]{Hou}
J.-M. Hou.
\newblock Light-induced {H}ofstadter's butterfly spectrum of ultracold atoms on
  the two-dimensional kagome lattice.
\newblock {\em CHN. Phys. Lett.}, 26(12), 2009.

\bibitem[HS84]{HeSjWsm1}
B.~Helffer and J.~Sj\"ostrand.
\newblock Multiple wells in the semi-classical limit {I}.
\newblock {\em Commun. in PDE}, 9(4):337--408, 1984.

\bibitem[HS87]{HeSjPise}
B.~Helffer and J.~Sj\"ostrand.
\newblock Effet tunnel pour l'\'equation de {S}chr\"odinger avec champ
  magn\'etique.
\newblock {\em Ann. Scuola Norm. Sup. Pisa}, Vol XIV(4):625--657, 1987.

\bibitem[HS88a]{HeSjHarper1}
B.~Helffer and J.~Sj\"ostrand.
\newblock Analyse semi-classique pour l'\'equation de {H}arper (avec
  application \`a l'{\'e}tude de {S}chr\"odinger avec champ magn\'etique).
\newblock {\em M{\'e}moire de la SMF}, No 34; Tome 116, Fasc. 4, 1988.

\bibitem[HS88b]{HeSjLecture}
B.~Helffer and J.~Sj\"ostrand.
\newblock {Equation de Schr{\"o}dinger avec champ magn{\'e}tique et
  {\'e}quation de Harper}, partie {I} champ magn{\'e}tique fort, partie {II}
  champ magn{\'e}tique faible, l'approximation de {P}eierls.
\newblock {\em Lecture notes in Physics (editors A. Jensen et H. Holden)},
  (345):118--198, 1988.

\bibitem[HS89]{HeSjHarper3}
B.~Helffer and J.~Sj\"ostrand.
\newblock Analyse semi-classique pour l'\'equation de {H}arper {III}.
\newblock {\em M{\'e}moire de la SMF}, No 39; Tome 117, Fasc. 4, 1989.

\bibitem[HS90]{HeSjHarper2}
B.~Helffer and J.~Sj\"ostrand.
\newblock Analyse semi-classique pour l'\'equation de {H}arper {II}
  (comportement semi-classique pr\`es d'un rationnel).
\newblock {\em M{\'e}moire de la SMF}, No 40; Tome 118, Fasc. 1, 1990.

\bibitem[JGT{\etalchar{+}}12]{JoETall}
Gy-B Jo, J.~Guzman, C.~K. Thomas, P.~Hosur, A.~Vishwanath, and D.~M.
  Stamper-Kurn.
\newblock Ultracold atoms in a tunable optical kagome lattice.
\newblock {\em Phys. Rev. Lett.}, 108(4):045305, 2012.

\bibitem[JZ03]{JZ}
D.~Jaksch and P.~Zoller.
\newblock {Creation of effective magnetic fields in optical lattices: the
  Hofstadter butterfly for cold neutral atoms}.
\newblock {\em New J. Phys.}, 5(56), 2003.

\bibitem[Ker92]{kerdelhue1}
P.~Kerdelhu\'e.
\newblock {\'E}quation de {S}chr{\"o}dinger magn\'etique p\'eriodique avec
  sym\'etrie d'ordre six.
\newblock {\em M{\'e}moire de la SMF, 2{\`e}me s\'erie, tome 51,}, pages
  1--139, 1992.

\bibitem[Ker95]{kerdelhue2}
P.~Kerdelhu\'e.
\newblock {\'E}quation de {S}chr{\"o}dinger magn\'etique p\'eriodique avec
  sym\'etrie d'ordre six : mesure du spectre {II}.
\newblock {\em Annales de l'IHP (section Physique th{\'e}orique)},
  62(2):181--209, 1995.

\bibitem[Mek03]{Mekata}
M.~Mekata.
\newblock Kagome: The story of the basketweave lattice. {Phys. Today Letters}.
\newblock
  \url{http://physicstoday.org/journals/doc/PHTOAD-ft/vol_56/iss_2/12_1.shtml},
  February 2003.

\bibitem[Mon13]{polyMontambaux}
G.~Montambaux.
\newblock {\em Conduction quantique et physique m{\'e}soscopique. Programme
  d'approfondissement Physique}.
\newblock {\'E}cole Polytechnique, avalaible in
  \url{http://users.lps.u-psud.fr/montambaux/polytechnique/PHY560B/PHY560B-2013.pdf}
  edition, 2013.

\bibitem[Nob]{nobel2010}
Nobelprize.org.
\newblock {The Nobel Prize in Physics 2010 - Advanced Information}.
\newblock
  \url{http://www.nobelprize.org/nobel_prizes/physics/laureates/2010/advanced.html}.

\bibitem[Nov11]{Novoselov}
K.~S. Novoselov.
\newblock {Nobel Lecture: Graphene: Materials in the Flatland}.
\newblock {\em Rev. Mod. Phys.}, 83(3):837--849, 2011.

\bibitem[Per60]{persson}
A.~Persson.
\newblock Bounds for the discrete part of the spectrum of a semi-bounded
  {S}chr{\"o}dinger operator.
\newblock {\em Math. Scand.}, 8:143--153, 1960.

\bibitem[PST06]{PSTd}
G.~Panati, H.~Spohn, and S.~Teufel.
\newblock Motion of electrons in adiabatically perturbed periodic structures.
\newblock {\em Analysis, Modeling and Simulation of Multiscale Problems}, pages
  595--617, 2006.

\bibitem[RS80]{ReedSimons4}
M.~Reed and B.~Simon.
\newblock {\em Methods of modern mathematical physics. IV. Analysis of
  Operators}.
\newblock Academic Press, 1980.

\bibitem[SBE{\etalchar{+}}04]{santos}
L.~Santos, M.A. Baranov, J.I. Cirac Jand H.-U. Everts, H.~Fehrmann, and
  M.~Lewenstein.
\newblock Atomic quatum gases in kagome lattices.
\newblock {\em Phys. Rev. Lett.}, 93(3), 2004.

\bibitem[Sim]{SimonXX}
B.~Simon.
\newblock Schr{\"o}dinger operators in the twenty-first century.
\newblock \url{www.math.caltech.edu/papers/bsimon/r40.ps}.

\bibitem[Sim83]{SiIHP}
B.~Simon.
\newblock Semi-classical analysis of low lying eigenvalues {I}.
\newblock {\em Ann. Inst. H. Poincar{\'e}}, 38:295--307, 1983.

\bibitem[Wil84]{wilkinson}
M.~Wilkinson.
\newblock Critical properties of electron eigenstates in incommensurate
  systems.
\newblock {\em Proc. Roy. Soc. Lond.}, A391:305--350, 1984.

\bibitem[Zak64]{zak}
J.~Zak.
\newblock Magnetic translation group.
\newblock {\em Physical Review}, 134(A1602-A1606), 1964.

\end{thebibliography}

\end{document}